\documentclass[a4paper,10pt]{article}
\usepackage[utf8]{inputenc}
\usepackage{amssymb,amsmath}
\begin{document}
\newcommand{\Lc}{\mathcal{L}_{c}}
\newcommand{\Rn}{\mathbf{R}^{n}}
\newcommand{\PSprod}{\breve{\otimes}_{\pi}}

\title{On uniformity for the polar to partially hypoelliptic operators}
\author{Tove Dahn}
\maketitle

\subsection{Terminology}
\subsubsection{Introduction}
 
 The main problem in this article, is to discuss the polar to hypoellipticy and to partial hypoellipticity
 relative $G_{8}$.
 Assume $Uf= \int f d U$, where $d U \in G$ and f is a symbol to a (pseudo-) differential operator.
 A necessary condition for hypoellipticity is absence of spirals in the polar set, for instance $U^{\bot} / U \rightarrow 0$ in $\infty$.
 Consider, for $d V \in G_{8}$ (\cite{Lie91}),  $K(G_{4})=(d V(x.y),d V(y,x))$, we denote $d V(y,x)= d V^{L}(x,y)$, . Assume $d W=x f_{x} + y f_{y}$, (that is $d W= d W^{L}$) then  $G_{4}$ is given by 
 $f_{x},x f_{x},y f_{x},x d W$ . Consider $G_{4} \rightarrow K(G_{4})=(G_{4},G^{L}_{4})$. With respect to
 conjugation $K(G) \rightarrow K^{\bot}(G)$, a necessary condition for hypoellipticity, is $K^{\bot} / K \rightarrow 0$
 or $K / K^{\bot} \rightarrow 0$.

 Assume $S_{j}(v)=(v, v^{L})$, $S_{j}^{\bot}(v)=(v,v^{\bot})$, relative different scalar products.
 Thus, when $S_{j}^{\bot} : A \rightarrow B$ and $S_{j} : C \rightarrow A$, with $\pi : C \rightarrow B$,
 it is necessary that we have existence of $(S_{j}^{\bot})^{-1}$ relative $S_{j}$.
 Assume $\int d S_{j}(v) = \int \rho d U_{1}$, where $d U_{1}$ denotes translation,
 we  note that  presence of lineality (sets of translation invariance) contradicts hypoellipticity. 
 $S_{j}^{\bot}$ is considered as a closed extension of $S_{j}$. We have that $S_{j}^{\bot}=S_{j}$ for large parameters, contradicts hypoellipticity.
  Given a regular analytic continuation, we have that $d S_{j}^{\bot} \rightarrow d S_{j} \rightarrow d V$
 is continuous.

 Consider $\Psi(t) : U(x,y) \rightarrow  {}^{t} U(x,y)$, so that $d {}^{t} U = \rho d U$.
 Assume $L_{1}$ a scalar product to $\Psi$ and $L_{2}$ a scalar product  to the orthogonal, assume $L_{1} \rightarrow L_{2}$
 algebraic. A closed set in the plane, according to $(U,U^{\bot}) \rightarrow (I,O)$ (cf. P-convexity),
 is through Tarski-Seidenberg mapped onto a semi-algebraic set in $(V,V^{L})f(x,y)=f(v,v^{L})$. Given $d I \notin G_{2}$, we do not have
 a closed (relatively closed) set in $(v,v^{L})$, in particular we do not have density. Note that given a strict condition,
 for instance $U^{\bot}$ algebraic, we have that $H(G_{2})$ generates $L^{1}(G_{2})$ (cf. \cite{Garding60}).  
 Thus, the problem is completely dependent of the topology, that is if $\mathcal{H}(G_{2})$ is generated by $H(K(G_{4}))$,
 then where we have symplecticity, $\mathcal{H}(G_{2}) \rightarrow H(K(G_{4}))$ is projective.
 We will us the notation $H(G)$ or $G_{H}$ dependent on the context.
 
 Assume $L^{1} \subset \mathcal{H}$. Given a set F is given by $\parallel \cdot \parallel_{1}=0$, it can be 
 considered as planar, that is generated by $G_{2}$. Example: assume $p$ a semi-norm, with $p(f) \leq C \parallel f \parallel_{1}$, that is we assume $f \in L^{1}$
 but $Vf \in \mathcal{H} \backslash L^{1}$, thus a neighborhood of $Vf-f$ is in $\mathcal{H} \backslash L^{1}$.
 Example: $f(v,v^{\bot}) \in H$ does not imply $f(v) \in H$. The polar set is given by $V^{\bot}f=0$, possibly
 non-removable.
 
 \subsubsection{The extended plane}
 
  Assume $d \varphi=dU$ non-decreasing and $f \in C^{0}$, we then have existence of 
 $\int f(x) d \varphi(x)=\int g(u)d u$, where g is undetermined over $u=const$, but does not affect action
 of U (\cite{Riesz27}).  
 
 \newtheorem{MP}{The moment problem}[section]
 \begin{MP}
 Assume $M=H(K(G_{4}))$ and $g \in C(G_{2})$.
 Assume, for a positive linear functional $B_{g}$, $B_{g}(g)=1$, for g fixed outside of $M$, we then have $B_{g}(f)=0$ for $f \in M$ implies $B_{g}(f)=0$ for $f \in C(G_{2})$. 
 Thus M approximates $C(G_{2})$ uniformly.
 \end{MP}
 (\cite{Riesz56})
 Assume $d V \in K(G_{4})$ and $f \in H(K(G_{4}))$, we then have $f(v) \rightarrow v$ continuous.
 Thus, given $g \in C$, we have that $f(v) \rightarrow v \rightarrow g(v)$ is continuous. When $K(G_{4})$ is dense 
 in X, we have that $g \circ v \in C$ iff $g \in C$. In particular $(K(G_{4}))_{L^{1}} \simeq (G_{2})_{L^{1}}$ according to Wiener.
 Consider $(K(G_{4}))_{H} \in d V \rightarrow \rho d V=d U$ a continuous continuation, with $\rho \in H$. Assume we have 
 $\sigma d V = d I$, with $\sigma \in H$. Consider the absolute continuous-convex closure $\mid d U(f) \mid \leq M \parallel d V(f) \parallel$
 where $d V  \in (K(G_{4}))_{H}$ is allowed to vary. Assume $g.l.b. \parallel d V(f) - g \parallel=d_{1}$, where $g \in C(G_{2})$, 
 Then $H(K(G_{4}))$ approximates $C(G_{2})$ uniformly iff $M \geq 1/d_{1}$. More precisely, assume $d U / \rho = \sigma d I$,
 where $\mid \sigma \mid \leq d$, we then have given $d U \geq d I$, $1/M \leq 1/\mid \rho \mid \leq \mid \sigma \mid \leq d$, that the inequality follows. The case $0 \leq d U \leq d I$ is obvious. Note that, given $\rho$ conformal, 
 we have existence of $1/\rho \in C$. Given
 $\rho \sigma$ reduced, we can represent $1/\sigma \rho$ in H, why $H(G_{2})$  approximates $C(G_{2})$.
 
 Consider $(K(G_{4}))_{H} \subseteq (G_{2})_{C} \in d I$, that is assume $f \in \mathcal{H}$ with continuation 
 $\tilde{f} \in C$, with the corresponding $g(v) \in H$ and $d V \in K(G_{4})$, that is $g(v) \rightarrow v$ continuous.
 Assume $H(G_{2})^{\bot} \simeq C(G_{2})$, in particular we assume the polar set has correspondingly $d V=0$ with $d V \in K(G_{4})$.
 Absence of non-trivial polar set for f, is motivated by density for $G_{2}$ regular.
 Example:  $f \in H(K(G_{4}))$,
 we then have $\int_{K(G_{4})} f dU=\int_{\Omega_{8}} f(u) d u=0$, when f is locally algebraic, $m \Omega_{8}=0$, (f is assumed complete).
 
 Reflection through
 $L_{j}$ defines uniformities according to K(d V) and $K^{\bot}(d V)$. Thus we have two diagonals
 $d K(V)=d K^{L}(V)$ and $d K(V)= d K^{\bot}(V)$. Given convexity for translates, we have that $G_{2}$ is decomposable to
 $K(G_{4})$. Given $d W=\rho(x,y) d I(x,y) \in K(G_{4})_{\mathcal{H}}$, where $\rho$ is topologically decomposable to
 $\rho=\vartheta (x) \otimes \vartheta(y)$, we have decomposability to $(G_{8})_{\mathcal{H}}$. Assume $\Delta=\{ d K(V)=d K^{\bot}(V) \}$ and $\Sigma=\{ d K(V)= d K^{L}(V) \}$, that is $\Delta \rightarrow \Sigma$ continuous,
 maps spirals on to symmetric sets. Obviously $K^{L}(d V)=K(d V^{L})$. 
 Let $\Psi(d U)=(d U, d {}^{t} U)$, that can be assumed harmonic.
 Note $\Delta U=0$ implies
 $d U \rightarrow d {}^{t}  U$ preserves type. Given $\Delta F=0$ on a domain, there is a corresponding G symmetric on
 the same domain. Example: $\phi(u,v)=c_{1} u + c_{2} v$ has $\Delta \phi=0$. 
 Starting from two scalar products, $L_{j}(f)=<g_{j},f>$, we have $L_{2}(L_{1}(f))=0$, for all $f \in (I)$,
 means either $L_{1}(f)=0$ or $L_{1}(f) \bot g_{2}$. $g_{j}$ defines measures regular over f.
 Further, $L_{1}(L_{2}(f))=0$, for all $f \in (I)$, gives that $g_{1} \bot g_{2}$ in measure space. 
 That is a symmetric neighborhood of $\Sigma$ $=\{ L_{2}=0 \}$ can define a polar set for transmission.
 In the same manner (geometric convexity), given $L_{1}(L_{2})=0$, a symmetric neighborhood of $\Delta$, can define a polar set for transmission.

 \subsubsection{A transmissions property}
 The transmission property (TP) that we discuss, is $(d V, d V^{L}) \rightarrow (d N , d N^{\bot}) \leftrightarrow (d U,d U^{\bot})$,
 where the middle element gives a rectifiable regular boundary, according to $(d N,d N^{\bot}) \rightarrow (dI,0)$. 
 A special case of TP is $T(\gamma_{1} \rightarrow \gamma_{2})(\zeta) \rightarrow \zeta$,
 that is the last element corresponds to a transformation in the $\zeta-$ domain. P-convexity for the group preserves hypoellipticity (\cite{Dahn2022}).
 Given the reversed last implication preserves topology, this transmissions property gives in particular
 an approximation property. Since our representations are based on the Stieltjes transform, instead of
 $C^{\infty}_{0}(\overline{X}) \rightarrow C^{\infty}(\overline{X})$ continuous (\cite{BdMv71}), we usually consider 
 $H(\overline{X}) \rightarrow C^{0}(\overline{X})$ continuous. Consider $T_{proj}(dU)=\{ \rho d U + \vartheta d U^{\bot} \}$.
 It is not obvious that $\rho \rightarrow 1/\rho$ preserves H, however if $\{ \rho =0 \}$ is of measure zero, 
 it is mapped on a set of measure zero, given the conditions above. Note in the case $H \in \rho \rightarrow 1/\rho \in C$,
 we have in $L^{1}$, that the order for the group is two, for the reversed last implication in TP.
 In order to avoid a spiral axis in the middle element, it is sufficient to consider $d N / d U^{\bot}=0$
 at the boundary. 
 
 Rectifiable boundary: (\cite{Riesz56}) Given $\Gamma$ with coefficients  u, v, w, we have $\Gamma$ rectifiable
 iff u,v,w BV (of bounded variation). For discontinuous curves with BV, we may have absence of two-sided limits (\cite{Garding77}),
 for instance $\lim d v \neq \lim {}^{t} d v$. Example: assume $I^{L} : (\xi,\eta) \rightarrow (\eta,\xi)$,
 $\overline{I} : (\xi,\eta) \rightarrow (\xi,-\eta)$, the harmonic conjugate satisfies $I^{\diamondsuit} \simeq \overline{I}^{L}$. 
 Thus given $d I^{L} U,d \overline{I} U$ correspond to closed forms and $d I^{L} U=\rho d I$, we have $\rho'$ can 
 be chosen as analytic. In particular, given $d I \notin G$, we may have $d U \in G$, but $d \overline{I} U \notin G$.
 
 Consider $T(\phi_{1} \rightarrow \phi_{2})$, where $\phi_{2}$ is analytic over a domain $\Omega$,
 so that $U \phi_{1}(x)=\phi_{1}(u) = \phi_{2}(x)$. 
 Thus we assume existence of $d U \in G_{2}$, so that $U \phi_{1} \in H$, 
 that is $\phi_{2}$ can be reached from every element  $\phi_{1} \in L^{1}$ through planar movements (outside the polar set). 
 Consider $dU=d V + d V^{L} \rightarrow (x,y,z)=\tilde{x}$. 
Assume $(\tilde{x}) \rightarrow (u_{1},u_{2})$ surjective. Thus
$f(u_{1},u_{2})=0$ implies existence of $\tilde{x}$, such that $f(\tilde{x})=0$. 
 
\newtheorem{tathet}[MP]{Density in the plane}
\begin{tathet}
When $d U=(d V, d V^{L})$ is dense in $d \tilde{x}$-space, 
that is $d U^{\bot}$ have trivial zero sets, 
we can reduce the two-mirror model to one scalar product. 
\end{tathet}

Assume $L_{j}(f)=<f,d U_{j}>$, with $d U_{j} \in G_{\mathcal{H}}$. 
In particular $L_{1}(f) \leq \Phi(f) \leq L_{2}(f)$ and given $L_{1} \rightarrow L_{2}$ projective,
the same holds for $L_{1} \rightarrow \Phi \rightarrow L_{2}$. Assume $d U^{\bot}=\rho d U$ and $\mathcal{H} \in f \rightarrow {}^{t} \rho f \in \mathcal{H}$
continuous, then
we have where $\rho \neq konst$, that $\Phi$ can be reduced to $L_{1}$, that is $\Phi(f)$ can be approximated
by $\int f d \Phi$, when $\rho \rightarrow 1$ regularly, even when $d \Phi \notin G$. 

Note that UI=IU does not imply $d I \in G$, even when $d U \in G$. 
Assume $GW$ has TP through completion. Example: $d I \notin G$, $GW \in d I$, $UW = I$, that is $U^{-1}=W$, with 
$d W \notin G$ and given $W \rightarrow I$ regular, we have completion to an approximation property, that is $UWI=IUW$.

Density is motivated by 
$f \in B_{pp}(K(v))$ (\cite{Schwartz57}). Thus, over the extended domain, $ \{ f(K(v),K^{\bot}(v)) < \lambda \} \subset \subset X$
implies $\{ f(u,u^{\bot}) < \lambda \} \subset \subset X$.
Given existence of an entire line L in $bd \Omega$
according to $f(K(v))=f(\tilde{x})$, we do not have hypoellipticity. P-convexity implies $K(d  V) \neq d I$ on any line L.

\subsubsection{Reduction to the boundary}

Obviously, $\mid f \mid = const$ does not imply $f = const$. Assume an approximation property,
$\mid \phi \mid=const$ has a sub-sequence, such that $\phi \rightarrow const$. 
This assumes $U \phi \rightarrow \delta_{0}(\phi)$
implies existence of V, such that $\phi(V \tilde{x}) \rightarrow \phi(0)$, that is $U \rightarrow I$ implies $V \rightarrow 0$.

A strict condition corresponds to $S : L^{1}(bd \Omega) \rightarrow \mathcal{O}(bd \Omega)$, (cf. \cite{BdMv85}) through
an infinitesimal movement, for instance S can be considered as a monotropic operator. Assume $S(P)=S(Q)$ with P,Q polynomial and
Q reduced at the boundary, then we have obviously (approximative) solvability, that is $S(P/Q)=I$. Thus
$f(u,v)$ corresponds to $P/Q d I(f)$ in $C^{0}$, that is by extending the domain, we can have solvability 
with respect to $P/Q$, without having it for P.
Assume $bd \Omega = \{ U=I \}$.
Given $bd \Omega$ is a subset of the domain for G, we have $d I \in G$. 

\newtheorem{ickenuclear}[MP]{Non-nuclear completion}
\begin{ickenuclear}
Assume $Uf - S_{1}(f)=f$, where $S_{1}$ is not nuclear, for instance a transformation that does not have a corresponding measure
$d V \in G_{H}$, that is $S_{1}$ defines a discontinuous representation.  Example: $S_{1}(f)-f \in C^{\infty}$ defines
a non-regular boundary over $S_{1}(f)=0$.  
\end{ickenuclear}

Example: the polar set is given by constant surfaces to f.
Example: the polar set is given by $(U-I)f^{N}=0$ implies $U^{\bot}f =0$.
Assume $rad (I) \simeq \overline{(I)}$, 
that is $(U-I)f^{N} \simeq (U-I)(f+c) = (U-I)f + c=0$ iff $(U-I)f=const$. 

Assume $d U = \rho d V$, with $d V=0$ on L and $\rho \neq 0$ on L. Note that given $\rho$ linear in V,
$d U = d V^{2}$ and given d V a projective operator, $d U \simeq d V$. Given all L are mapped on boundary cycles,
we have that projective measures preserve the boundary. 

Assume $\Omega$ semi-algebraic and $\Omega \rightarrow bd \Omega$ continuous. Example: f regular and $f \neq 0$,
 $0 = \int f(u)d u=\int f(x,y)d U(x,y)$
with $(x,y) \in L$, for instance $u \sim 0$ (homotopic with 0), that is we consider $If = f(0)$. Given L
is defined by $U=V$, we have $f(u,v)-I$ is regular. When $U$ is a spiral, the mapping $\Omega \rightarrow bd \Omega$
is not continuous.

Nuclear topology for $\mathcal{H}$ and $d U \in G_{\mathcal{H}}$, implies $UI=IU$, in particular $ U-I$ is algebraic,
that is $(U-I)I=I(U-I)$ and given $d (U-I)=P dI$, for P polynomial, through Hurwitz, L must be trivial.

Assume $d V(f)=f_{x}$ and $d V^{L}(f)=f_{y}$, given $f$ symmetric in x,y, we have 
$f_{x}=0$ iff $f_{y}=0$, that is $f(v,v^{L})$ symmetric. Note $d (V,V^{L}) \in K(G_{4})$
iff $d (V,{V^{L}})^{\diamondsuit} \in K(G_{4})$ (harmonic conjugation). Example: $(f_{x})^{\bot}=-y f_{x}$ and $(f_{y})^{\bot}=x f_{y}$
that is $(f_{x} + f_{y})^{\bot}=-y f_{x} + x f_{y}$. 

Assume the boundary is defined by a normal operator $(d N,d N^{\bot})$, where $d N \rightarrow d I$ iff $d {}^{t} N \rightarrow d I$,
that is we consider a symmetric neighborhood of $d I=d {}^{t} I$. Further, $d N^{\bot} \rightarrow 0$ regularly
at the boundary. Consider the continuation to $(d U,d U^{\bot})$, why $d (UN)^{\bot} \rightarrow 0$ regular.
A two-sided regular limit, assumes $\lim UN = \lim U {}^{t} N$. Consider now $d V \times d V^{L} \rightarrow d N$,
so that $d V \rightarrow d I(x)$ iff $d V^{L} \rightarrow d I(y)$, which defines a symmetric
neighborhood in (x,y) and given that translates are dense in the domain $d V + d V^{L} \rightarrow d I(x,y)$ (\cite{Garding77}).
\subsubsection{TP relative 1-parameter group}

Consider $\gamma(\zeta_{T})=\gamma_{T}(\zeta)$, as a topological completion.
Assume $d U \rightarrow d V \rightarrow d V_{1}=\rho_{1} d I=0$ regular, that is dU has an analytic 
base through topological completion, where through translation $\rho_{1} \neq 0$ locally.
Concerning projective mappings in this case, $d V_{1}(\gamma_{1})=0$ in $\mathcal{D}_{L^{1}}$, implies existence of $d V_{2}$, 
so that $d V_{2}(\gamma_{2})=0$ in H, when $\gamma_{j}$ are entire lines.

Assume $W f_{0} - f_{0} + f_{0} - I \sim 0$,
then $W$ is very regular, that is $W f_{0} - f_{0} \sim 0$, given $f_{0}$ very regular. Note that given strict pseudo convexity, 
that is given v locally algebraic and $f(u,v) \rightarrow  f(0)$, 
we have a transversal approximation. 

Assume $d U^{\bot} = \rho d U$, where
$\rho \equiv const$ gives the diagonal. In particular, $d^{2} U^{\bot} = \rho' d U + \rho d^{2} U$, 
that is given $\rho \neq 0$, we have that $d U^{\bot} \rightarrow d U$ preserves corresponding closed forms (projective mapping)
iff $\rho' d U=0$, thus projective mappings include spirals. 

In analogy with the fundamental theorem in algebra:
$f(u,v)=0$ modulo $C^{\infty}$, means that given f algebraic with respect to (u,v), we have convergence in $C^{\infty}$.
Assume UV=VU defines $Z=\{ (u,v) \quad f(u,v)=0 \quad u \in G \quad \exists v \in G \}$
a domain for convergence in $C^{\infty}$. Polynomials approximate $C^{\infty}_{c}$ (\cite{Treves67}, Cor. 2 Lemma 15.1), that is
given $U,V$ pseudo local and given $f \in C^{\infty}_{c}$, we have existence of $P(u,v)=0$ modulo $C^{\infty}$ over Z
and $\mid f - P \mid < \epsilon$.

 \subsubsection{A weaker strict condition}
Consider $X' \subset X$, for instance $X'$ polynomial, that is $v \in X'$
does not imply $1/v \in X'$. Given $v=y/x \rightarrow 0$ or $x/y \rightarrow 0$, we have 
$1/(v+1/v) \rightarrow 0$ in $\infty$. Thus we can consider $X' \times X \in (v,1/v)$ as one sided, 
that is given density for translates, we can consider $v + 1/v$ as approximately in $X'$. 

Example: Assume $\Sigma=(\frac{U}{V} + \frac{V}{U})$, with $\Sigma d \Sigma=d I$ and $d \Sigma(x,y)=d \sigma$.
We then have, given $u/v \rightarrow 0$ or $v/ u \rightarrow 0$, that $1/\sigma \rightarrow 0$
and $d \Sigma^{2} \sim_{0} d I$ (geometric equivalence).
Given U,V projective operators, we thus have that $(V^{-1}U + U^{-1} V) F(x,y,z) = F(x,y,z)$ implies $U+V=UV$.
Given $d U \rightarrow d V$ projective, we have weakly $U \rightarrow V$ projective. Given d U 
sub-nuclear, we have for instance $d U \rightarrow d U^{2}$ projective. Thus, $d U \rightarrow d (I-U)^{-1}$
projective.

Given $d V^{-1}U \simeq (p/q) d I$, we have $pq / (p^{2} + q^{2}) \rightarrow 0$.
The condition p/q $\rightarrow 0$ or q/p $\rightarrow 0$, is sufficient for 
$L_{p^{2} + q^{2}} \subset L_{pq}$. 

\newtheorem{HC(G)}[MP]{Preservation of convexity for weak Lie type}
\begin{HC(G)}
 Problem: consider d U, such that the type is preserved for the convex closure of the domain, 
 in particular in the case where $d I \notin G$. 
\end{HC(G)}

Consider $T_{ac}(G_{H}))=\cup_{\mid \rho \mid \leq 1} \rho d U$, where $d U \in G$ and $f \rightarrow \rho f$
continuous. $C(G_{2})$ allows cluster-sets. Example: $\mid f \mid =1$ allows infinitely many values for
$U_{2}$. A complete limit is dependent of $U_{2}$.

Define $\int d <d E,f>=\int d E(f) = E(f)$, that is we assume $<d E,f>$ in $L_{ac}$.
Given $d <dE,f> = \rho <d I,f>$, we have that where $\rho \neq konst$, we can not conclude $d E = d I$. Further, constant surfaces give the type, but $<d E,f> =<d I,f>$ does not imply 
$d E = d I$ in regular topology, for instance E can have a non-trivial kernel. 

Given the convex closure in $C(G_{2})$, it is not obvious that the type is preserved. Example: $Uf \rightarrow f$ iff
$f(u) \rightarrow 0$ or weaker $f(u+v) \rightarrow 0$, with $u / v \rightarrow 0$ in $\infty$. Assume
$d U = \rho d V= d WV$, then in order for $U=I$ to imply $V=I$, we have to determine W, so that $WV=VW$. Example: $\rho$ bounded and absolute continuous,
we then have $d \rho=0$ iff $d (1/\rho)=0$, that is $\rho=const$ iff $1/\rho=const$. In particular $\rho d V \rightarrow dV$ is
absolute continuous, further projective.

Assume $(U,V)$ a continuation of projective movements through conjugation. Given U,V of the same type,  
for instance $d (U+V)=\rho d U$ and U+V=I, the conjugation gives a maximal domain 
for invariance. Assume U,V of different types with U+V=I. Assume
$Uf = Vf$ implies f=0. Assume $\tilde{U}f$ continuation with 0 and in the same manner for $\tilde{V}$, 
we then have $\tilde{U} + \tilde{V}=I$, with $\tilde{U}f=\tilde{V} f$ implies $f=0$, that is given U,V with one-sided support,
U,V can be continued to a projective conjugation. 

Consider
$d V(\rho f) \simeq {}^{t} \rho d V(f)$, given $\rho \rightarrow {}^{t} \rho$ preserves $T_{ac}$ - convexity
we have that dV is defined on the  $T_{ac}$-closure, for instance this is the case when $\rho$ polynomial. 
However, $\rho \rightarrow 1/\rho$ does not necessarily preserve
convexity, a sufficient condition is $\rho$ absolute continuous. Note that $(I) \subset (J)$ does not necessarily 
imply  $T_{ac}(I) \subset T_{ac}(J)$.
A necessary condition for inclusion, is that the corresponding coefficients have 
$ \rho_{J}/ \rho_{I} \rightarrow 0$ in $\infty$

In the plane, we have that $W + W_{1} + R$ according to TP, preserves convexity (modulo removable sets). 
Assume W algebraic, $W_{1}$ preserves one-sidedness, for instance reflection and R negligible.
Consider the $T_{ac}$-closure, $T_{ac}(I)=$
$\{ \rho f \}$ of $f \in (I)$. In the case $d I \notin G_{2}$, we consider $d V \rightarrow d V \times d V^{L} \rightarrow G_{2}$
where $d V \in G_{8}$ is the result of the inverse scheme, that is $d U \rightarrow d V$ is two-valued. 
Assume the condition one-sidedness in $(d U)$ is mapped onto one-sidedness in (dV). Then TP can be 
continued to $G_{8}$.

Assume d (U+V) $\rightarrow d I$, with dV compact ($C^{\infty}$), meaning the integral operator defined by dV is compact.
Then, in the same manner, d U - d I completely continuous
corresponds to dU Fredholm. $d (U-V)^{2}=d U^{2} + d V^{2} - 2 d UV$, with $d UV \rightarrow d I$,
thus $2 U d U + 2 V d V - 2 d I \sim d(U-V)^{2} \rightarrow 0$, that is $U d U + V d V \rightarrow d I$.
In particular $(U-V) d (U-V) \rightarrow 0$, where the diagonal gives non-regular approximations. 

Example: assume $d U \rightarrow d W$ are conjugated, according to $d U(f)=-\{ g,f \}$ and $d W(f)=-\{ h, f \}$. 
Assume $\Psi(g)=h$, given $\Psi$ preserves linearity and order, we have
that the type is preserved. In particular starting from $d U \rightarrow d {}^{t} U$, we have 
$UT(\phi)=T(U \phi) + (-\Delta g)T(\phi)$, that is the type is preserved under transponation, given $\Psi$ preserves
harmonicity. 

\subsubsection{Inverse lifting principle}

\newtheorem{vectorvalued}[MP]{Vector valued distributions}
\begin{vectorvalued}
Consider $d U \rightarrow d V$ conjugation and $\Phi(d U) \rightarrow (d U,d V)$ and
for a second conjugation $d U_{1} \rightarrow d U_{2}$, $\Phi^{2}(d U)=(d U_{1},d U_{2})$. 
With this construction, $\Phi^{2}$ can be algebraic,
without $\Phi$ being algebraic. Assume in particular, $\Phi_{1} : d W \rightarrow d {}^{t} W$, $\Phi_{3}: d W \rightarrow d W^{\bot}$
and $\Phi_{2} : d W \rightarrow d W^{*}$ (algebraic dual). Using (\cite{Schwartz57}) $<<T,\phi>_{1} , d U>_{2}=<<T,d U>_{2},\phi >_{1}$ 
we have $< T,d W^{*}>_{2}(\phi)=<T(\phi), d W>_{2}$. Consider completion to $d \tilde{W}^{*}$ closed,
so that $d \tilde{W}^{*} \rightarrow d W^{\bot}$ linear (single valued), that is $d W^{\bot}=\rho d \tilde{W}^{*} = \rho \sigma d W^{*}$,
or $d W^{\bot}=d \tilde{W}^{* k}$.
\end{vectorvalued}
 
Assume $\sigma$ regular, that is $d W^{\bot}=0$ iff $d W^{* k}=0$, we then have $d W^{\bot}=\sigma^{k} d W^{*}$,
given $\sigma^{k}$ algebraic, the continuation has removable singularities in a regular neighborhood. 

Consider $\mid \Phi(f) \mid^{2}=(\mid d V(f) \mid, \mid d V^{L}(f) \mid)$, where 
$\int \mid \Phi(f) \mid^{2}$ finite, defines a normal surface, (\cite{Nishino62}), that is 
an approximations principle for $\mid \Phi(f) \mid^{2}$, does not assume $\Phi$ algebraic. 
Consider $d \sigma \sim \mid X \mid^{2} + \mid Y \mid^{2}$. 
 Consider $d \sigma=(1/v) d xdy$, where $X,Y$ are algebraic, that is $1/v=0$ removable.
 Starting from $-\xi Y + \eta X$, we have that v determines the character of movement. 

Assume $d U = d V + d V^{L}$ and $d U = d I(x,y)=d I (x) \otimes d I (y)$. We then have $dV=d I(x)$
and $d V^{L}=d I(y)$ given decomposability. That is, $d V=0$ implies $d V^{L}=0$, 
that is we assume $d U$ preserves zero-lines in the sense of $x \rightarrow y$ projective. 
Note, that given convexity and a strict condition
analogous with $B_{pp}$, we can assume $(v + v^{L})$ generates $(v, v^{L})$.

 \newtheorem{MPIII}[MP]{MP conjugation}
 \begin{MPIII}
 Assume for $\rho$ not necessarily polynomial, $q_{1} \leq (\rho + 1/\rho) \leq q_{2}$, where $q_{1},q_{2}$ polynomial.
 Then $d U_{\rho} + d U_{1/\rho}$ has an approximation property.
 \end{MPIII}
 
 Obviously $q_{1} d I \leq (\rho + 1/\rho) d I \leq q_{2} d I$, that is $\int_{\Omega} q_{2} d I=0$ implies
 $\Omega$ a zero-set. Assume $d U \rightarrow d U^{\bot}$ according to $d U + d U^{\bot}=(\rho + 1/\rho) d I$,
 that is as long as $d U + d U^{\bot}$ BV, we have convergence for d I in $\infty$. 
 
 Note that, given $\rho=e^{\phi}$, we have $\rho + 1/\rho \sim 2 \cos \phi$, that is a harmonic representation.

\subsubsection{Very regular approximations}
We assume the condition $d I \notin G$ corresponds to $0 \notin $ the domain for $T \in \mathcal{H}$, 
that is an annular region. 
Consider $f(u,v) \simeq f_{1}(u-v)$, when (u,v) is not tangent to the diagonal, 
we do not have a vector field (\cite{Malgrange71}).
We assume $d I \notin G$ defines discontinuous
approximations.

\newtheorem{extended}[MP]{Regularity in an extended plane}
\begin{extended}
Example: $d I \notin (G_{2})_{C^{\infty}}(u_{1},u_{2})$, but $d I \in (G_{3})_{C^{\infty}}(u_{1},u_{2}, du_{3})$.
$f(u,v) \notin C^{\infty}$,
but we have existence of W, such that $f(u,v,w) \in C^{\infty}$. Example: f very regular, that is 
$f(u,v) - I(u,v) \in C^{\infty}$ iff $f(u,v,w) \in C^{\infty}$, that is $Wf=0$ (modulo $C^{\infty}$) over u=v.
\end{extended}

In the last example, W corresponds to, $d W \in G_{reg}$, that is regularization in a complete limit.
Note that $f(u,v) \rightarrow f(0)$, given $d I \in G(u,v)$. Given f algebraic, when $U \rightarrow I$,
f must be algebraic in x,y,z, in this case.

\newtheorem{Hausdorff}[MP]{Uniformity}
\begin{Hausdorff}
Assume U,V sub-nuclear (modulo $C^{\infty}$) and projectively conjugated, so that $U \rightarrow I$ is regular iff $V \rightarrow 0$ regular.
We then have
$E- I \in C^{\infty}(x,y,z)$ implies $E-I \in C^{\infty}(u,v)$. Conversely $f \in C^{\infty}(u,v)$,
implies $f \in C^{\infty}(x,y,z)$, where $d U(x,y,z)$ is analytic.
\end{Hausdorff}

The proof assumes a separation property for d I.
It is sufficient for this to use sub-nuclear (pseudo local) movements. d I is sub-nuclear since $C^{\infty}$ is nuclear.
Consider $\int f d U=\int f(u) d u=Uf$, it is sufficient that $d U \in C^{\infty}$ in x,y,z (analytic).
Consider $\pi (x,y,z)=(u,v)$ and $d \pi=(d u,d v)$. Example: $d U(x,y,z)=\mu d x d y d z$, with $\mu \in C^{\infty}$.
Example: $(x,y,z) \rightarrow u \rightarrow (u,v)$, that is dv regular measure in nbhd R(d U) (the range of dU), 
corresponds to a separation axiom and d u=d v corresponds to a boundary for regular continuation.

\newtheorem{TPex}[MP]{TP and preservation of convexity}
\begin{TPex}
 We can determine d W, so that $F(w)=F(u,v) - I(u,v)=0$ (modulo $C^{\infty}$). Continue W to $\tilde{W}F=P(U,V) F - Q(U,V)I$ compact and regular
(P,Q linear in u,v), that is $\tilde{W}F=0$ iff $\frac{P}{Q}F(u,v)=const I(u,v)$, that is P/Q preserves weak constant surfaces iff
$\tilde{W}F=0$. 
 
\end{TPex}

Example: assume $d W(F)(u,v)=d I(u,v)$ (modulo $C^{\infty}$), that is given F solution to a hypoelliptic operator d W, 
we have existence of W with 
$F(w)=F(u,v)-I(u,v)$.
Consider $\Phi : \mathcal{H} \rightarrow L^{1}$, defines a topology that can be reduced to $L^{1}$
through two reflections. Consider $\rho d U + \vartheta d V \in \tilde{G}_{2}$,
where we choose $\rho, \vartheta$, so that $\tilde{G}$ approximates G uniformly.
Assume $P(U,V)f=g$, when $U=I$, $P(V)f=f(P(v))$, given $f(0)=0$
we have $P(v)=0$ implies $g=0$, given $f \neq 0$ close to 0, g has algebraic zeros. 

Note that $d I(u,v)$ with support on the diagonal, gives an approximation property with cluster sets on
the diagonal (a ``spiral''). With this representation the fundamental solution may not be representable in the plane.
For f phe and the corresponding parametrix E, we have that f is hypoelliptic outside the diagonal, 
also on the kernel to E (\cite{Dahn13},\cite{Dahn18}).
This implies uniformity for $rad (I_{HE})$. For more general f, the diagonal is not necessarily 
the only counter example to hypoellipticity.

\subsection{Fundamental concepts}
\subsubsection{Convexity}
 
 For holomorphic convexity, we consider $<g,f> \rightarrow G(f)$, where G is holomorphic, but not necessarily linear or nuclear,
 that is we do not assume an approximation property. Example: $G(f)=e^{P/Q}(f)$, with for instance Q reduced, 
 which implies TP (\cite{Dahn15}). Given $d I \in (G_{2})_{C}$,
 simultaneously $d I \in (G_{2})_{L^{1} \cap C}$, it is sufficient to consider $\parallel f \parallel=1$.
 Finite dimensional normable spaces are nuclear (\cite{Treves67} Prop. 50.2 Cor. 2).
 
 Assume T the symbol for a parametrix, and $T-I=e^{E}$. Sufficient for T-I to have negative
 type, is $e^{\mid x \mid} \mid T-I \mid \leq C$, when $\mid x \mid \rightarrow \infty$. Assume $S(d I f)=x f_{x} + y f_{y}$.
 Note that, $\mid S(d I) \mid \leq \mid x \mid \mid d I \mid + \mid y \mid \mid d I \mid$.
 $E( d V(\phi))=<E, d {}^{t} V>(\phi)={}^{t} VE(\phi)$ and $\mid <E, S(d I)>(\phi) \mid \leq \mid E(\phi) \mid (\mid x \mid \mid d I \mid+ \mid y \mid \mid d I \mid)$.
 Denote $\tilde{V} (T-I) =e^{V E}$, given V algebraic, we have
 $\tilde{V} \sim V$. Given $\mid x \mid \mid E(d I) \mid \leq C$ and $\mid y \mid \mid E(d I) \mid \leq C´$,
 we have $E(d I) \rightarrow 0$, that is
 corresponding to type 0. Thus, the kernel to $T-I$ is trivial in $\infty$.
 Conclusion: given T-I has linear representation in the phase, 
 $T-I = e^{E}$, where (T-I) has type 0 and given $E(S(d V) (\phi))$ is locally bounded, when $d V \rightarrow d I$ regularly, 
 we have that $\tilde{V} (T-I)$ has negative type.
 Negative type gives a nuclear topology.

\newtheorem{holje}[MP]{Geometric convex closure}
\begin{holje}
Consider $<f , d U> \leq C_{G}<f, d I>$, $\forall d U \in G$ and a constant $C_{G}$,
dependent on $f \in (I) \subset C^{\infty}$,
that is we assume d U BV, but not necessarily $dI \in G$.
\end{holje}

Example: $d U = \rho d I$, where $\mid \rho(x,y) \mid \leq C_{G}$ on a compact set.
The sets $X=\{ f \} \subset \subset C$ in the domain for G, are geometrically convex iff $\tilde{X} \subset \subset C$. 
Note that $f \equiv 0$ over $G_{8}$, implies $f \equiv 0$ over $G_{2}$.
We start with two scalar products $S(v),S^{\bot}(u)$, and $G_{8} \rightarrow G_{2}$ : $S^{\bot} S(v)=
((v,v^{L}),(v,v^{L})^{\bot})$.

\newtheorem{relativ}[MP]{Relative convexity}
\begin{relativ}
Assume $\rho$ analytic,
with absence of essential singularities in $\infty$, then there are $d U,d V$ analytic, so that $d U = \rho d V$. 
In particular,
given $\frac{\delta}{\delta v^{\bot}} \rho(v,v^{\bot})=0$ and $\frac{\delta \rho}{\delta v} \geq 0$, we have that $U$ is relatively convex. 
Given a strict condition,
we can continue to $\rho \in L^{1}$. 
\end{relativ}

Let for instance $d V^{\bot}=\sigma d (V-I)$.
Assume $d U = \rho d V$, where $\rho'(v) \geq 0$ locally, that is $\rho$
monotonous and single valued, thus in the plane linear. Assume $d U_{2}=\mu d U_{1}$, with $\mu$ analytic and $\rho d U_{1} + \vartheta d U_{2} =d I$,
that is $(\rho + \vartheta \mu) d U_{1}=d I$. In particular for $1/(\rho + \mu \vartheta) \in (I)$, a multiplier ideal,
we can discuss $\rho + \vartheta \mu=const$. Example: determine $\rho, \vartheta$, 
so that $\overline{\delta} \rho + \vartheta \overline{\delta} \mu + \mu \overline{\delta} \vartheta=0$,
for instance $\overline{\delta} \rho / \overline{\delta} \vartheta=-\mu \neq 0$.
Given $\rho d U_{1} \rightarrow d I$ iff $\vartheta d U_{2} \rightarrow 0$, we get a generalization of 
projectivity. Example: $\vartheta \rho (\frac{1}{\rho} d U_{2} + \frac{1}{\vartheta} d U_{1})=\rho d U_{1} + \vartheta d U_{2}$
for instance $\vartheta=1/\rho$. Note that $\frac{1}{\rho} + \frac{1}{\vartheta}=1$ iff $\rho \vartheta=\rho + \vartheta$.

Concerning $G_{8}$, given for instance $d V(f)=f_{x}$, holomorphic
convexity is not preserved, consequently not geometric convexity.

\newtheorem{linjarisering}[MP]{Linearization}
\begin{linjarisering}
Assume d U
(projective), completed through dW to UW algebraic (analytic), then type is not necessarily preserved.
Thus, we can have an approximation property for $T_{ac}(d U)$, that is not present for dU.
\end{linjarisering}
Example: $d UW=\rho d U$, with $\rho \neq 1$, where $d U = d I$.
.
Example: irreducibels:  $d U = \rho d W=d VW$, we then have $d U=d I$ implies $d V=I$
or $d W=d I$, that is $IW=VI$, corresponding to continuation
under preservation of type.

Starting with two scalar products, assume $UW$ closed,
where $W^{\bot}$ is defined relative $S_{2}$, so that $U^{\bot}W^{\bot}=W^{\bot} U^{\bot}$.
Thus, we can define $(UW)^{\bot}$ through annihilators, that is not relative type of movements.
Alternatively, consider completion to a closed movement $(d V,d V^{L})$, 
further $((),()^{\bot})$ relative normed spaces and annihilators.

Assume convergence according to  $\mid Uf \mid \rightarrow \mid f \mid$, this does not imply $\mid Uf - f \mid \rightarrow 0$,
but given existence of d V, so that $\mid VU f - f \mid \rightarrow 0$, we have mean convergence, however
dependent of $(U,V)$. Example: assume $S=UT -TU$, then
a conjugated movement relative mean convergence, is given by $VS=0$ on a given set $\Sigma$. 
A polar set can be given by $VS=0$ implies $d V \notin G$.

Example: Assume $d U=\rho d V$, with $\rho \in C^{\infty}$, then we have given monotropy (isolated singularities), 
existence of $v$ analytic close to $\rho=1$.
Assume $UV \rightarrow I$ regular with $d V=(1/\mu) d I$ and $\rho / \mu = v$ analytic, that is it is not necessary
to assume $1/\mu$ analytic. The type for dV is given by $1/\mu=const$.

 \subsubsection{Symmetry}
\newtheorem{symmetri}[MP]{Symmetry mapping}
\begin{symmetri}
Assume $\sigma$ is defined relative a scalar product, for instance $U \rightarrow {}^{t} U$
and $\sigma : (V, V^{L}) - {}^{t} (V,V^{L})$.
Symmetry implies
$\sigma \sim 0$, for instance homotopic with zero. Further we assume $\sigma$ is a projective
mapping, that is maps zero lines onto zero lines. 
\end{symmetri}

Note that $U \rightarrow {}^{t} U$ does not necessarily preserve type, that is $U \rightarrow I$
does not imply ${}^{t} U \rightarrow I$, that is $\Sigma(U) \neq \Sigma({}^{t} U)$
regular. 
Consider $\sigma_{1}$ with respect to $d V \rightarrow d V^{L}$ and $\sigma_{2}$ with respect to  $d V + d V^{L}=d U \rightarrow d {}^{t} U$.
For symmetric functions, $\sigma_{1}$ does not define a regular symmetry. Obviously, when $\sigma$ symmetric, $\sigma_{1}(d V) F(\phi) \simeq F (\sigma_{1}(d V) \phi)$.
Given F linear and $U,{}^{t} U$ sub-nuclear,
we have $F (\sigma_{2}(d U) \phi) \simeq \sigma_{2}(d U) F(\phi)$.

Projectivity can be compared with preservation of closed forms, these are locally exact. Given $d U,d U^{\diamondsuit}$ correspond to closed forms,
using the involution condition,
we have analyticity for d U (\cite{AhlforsSario60}). Assume $d U = \rho d I$ and $d U^{\diamondsuit}=\vartheta d I$. 
On $\rho \neq 0$ and simultaneously $\vartheta \neq 0$, projectivity is preserved.

\newtheorem{twomirror}[MP]{The two-mirror model} \label{two}
\begin{twomirror}
TP : $A \rightarrow_{v} B \rightarrow_{\pi} C$,
where $v A \simeq \pi^{-1} C$, that is given $\pi$ invertible, we have a lifting principle. 
A strict condition is sufficient for $(u^{\bot}) \rightarrow (u) \rightarrow ( x )$ continuous.
Given an approximation property, we have $(u) \rightarrow (x)$ regular.
\end{twomirror}

Assume f bounded and analytic, symmetric on a planar domain with a simple boundary L, then f is symmetric on a disk (\cite{Collingwood66}). 
Example:
$E(u,u^{\bot}) - I \sim 0$ on $L = \{ u=u^{\bot} \}$, then E is not analytic.
Further, given $(U,U^{\bot})$ discontinuous, then $E(u,u^{\bot})$ is not analytic.

 \subsubsection{Discontinuous approximations}

 Assume $f d I$ BV and analytic, where $f=0$ on $L \subset bd \Omega$ implies $f \equiv 0 $ on one side of L.
 Consider $\pi : (x,y,z) \rightarrow (u,v)$ and $L=\{ u=v \}$, given $f =0$ (modulo $C^{\infty}$) on L, 
 implies that $f \in C^{\infty}$, for instance when $v/u \rightarrow 0$. The condition $fd I - dI =0$ on L leaves space for a spiral.

Note that $f$ one-sided in $G_{2}$, does not imply that
$f$ is one-sided in $G_{8}$, but one-sidedness in $K(G_{4})$ implies one-sidedness in $G_{2}$.

Assume f bounded and =0 on $L \subset bd \Omega$ implies $f \equiv 0$ on a half space (convex set), 
that is given I gives reflection of supp f through $bd \Omega$, $If \cap fI = \{ 0 \}$, 
that is when If defines linear continuation,
we have orthogonal decomposition of the corresponding measure. With these conditions, projective movements have a decomposition $U_{+} - U_{-}=U$.
Assume f separately linear in u,v, we then have $<d V, <d U,f>>=<d U \times d V,f>$, that is f linear in $d U \times d V$.
Note that $(U+V)f=f(u+v)=f(u)+f(v)$, only when f linear. 

\newtheorem{cluster}[MP]{Linear phase}
\begin{cluster}
On a geometrically convex set, we have existence of $C_{K}$, so that
$L_{1}(L_{2}(f)) \leq C_{K}$ iff $L_{2}(L_{1}(f)) \leq C_{K}$. 
\end{cluster}

Assume $e^{L_{1}}$ defines
an algebraic symbol, why $\{ e^{L_{1}} \leq C \}$ corresponds to a semi-algebraic set. Note that
$L_{1} \equiv 0$ (or constant) on an entire line, gives cluster sets,
that is we do not have relatively compact sets. 

Assume f analytic, with $d U(f) + d U^{\bot}(f)=0$, further $d U \bot d U^{\bot}$.
Consider $F(u,v) \rightarrow u$. Assume $\Sigma(U)=\{ U - I =0 \}$. To determine type of movement,
it is necessary that $\Sigma \cap$ the domain for F. Assume $U \rightarrow I$ implies $V \rightarrow 0$, 
then we have $(U+V)F \rightarrow F$ on $\Sigma$. Thus, when the conjugation is
given by $\Sigma(U)=N(V)$, the movements can be identified regularly.

\subsubsection{The stability group}

A stability group is formed by $g_{0}=Id$ and $g_{i} g_{0}= g_{0} g_{i}$, where $g_{i}$ is a connected
family, not necessarily a subgroup. (\cite{Cartan41})

When the stability group gives two-sided regular approximations, we have that $g_{i} I = I g_{i}$ implies $g_{i}$ (topologically) algebraic,
Example: $F(\gamma) \in L^{2}$, $F(\gamma \rightarrow \gamma') \in \mathcal{D}_{L^{2}}$ and
$F(\gamma \rightarrow \gamma' \rightarrow \gamma'') \in H$, gives a strict condition, a completion
to stability.

A subgroup approximates 0, that is stability defines a (topologically) algebraic subclass, that is not necessarily
an analytic subgroup.  Example: UI=IU implies (U-I)I=I(U-I),
where $d U - d I$ can define a subgroup, even when $d U$ does not.

Two-sided limits can be motivated by regular complements, for instance $d U$ analytic with existence of d V analytic
and $d U \bot d V$. Assume $d U = \rho d V$, with $\rho \neq 0$ conformal in the plane,
that is starting from $d V \rightarrow d V \times d V \rightarrow d U \times d U^{\bot}$ regular and
conformal, that $d U \bot d U^{\bot}$, for some $d U^{\bot}$.
Note that conformal mappings are locally invertible, that is the last mapping in the two mirror model (section \ref{two}) is reversible.

\subsubsection{Uniformity for the polar set}

\newtheorem{dummy}{Uniformity}[section]
 \begin{dummy}
 Assume $E(f)=I(f)$ (modulo regularizing action) localization of the symbol, then E(f)=0 represents the polar set.
 $E(f)=\int f d E$, with $d E - d I=0$ (modulo regularizing action) over $E(f)=0$, implies that the diagonal
 uniquely gives the polar for regularizing action. Thus the representation is Hausdorff uniform (\cite{Bourbaki89}).
 \end{dummy} 

 Note that E corresponding to the parametrix to a partially hypoelliptic differential operator, satisfies the conditions for 
 Hausdorff uniformity relative $K(G_{4})$ above (\cite{Dahn18}).
 Concerning $G_{8}$: the six first are not complete and thus do not preserve for instance uniformity, compactness. 
 For $d V_{1}=x^{2} f_{x} + xy f_{y}$, $d V_{2}=xy f_{x} + y^{2} f_{y}$,
 we have $d V_{1}/d V_{2}= \frac{\frac{x}{y} f_{x} + f_{y}}{f_{x} + \frac{x}{y} f_{y}}$.  Example: $(x^{2},xy) \rightarrow (x^{2},y/x)$
 and $(xy,y^{2}) \rightarrow (y^{2},x/y)$. $F(x,y) \rightarrow F(x^{2} + y^{2}) M(y/x,x/y)$,
 with M algebraic. This representation defines F uniformly, where M=const.

 Note $d U^{\diamondsuit} \bot d W^{\diamondsuit}$, implies $d U \bot d W$. For an orthogonal relation,
 we assume the coefficients to the movements are related through a scalar product, for instance Legendre. 
 We do not necessarily assume that $d U^{\bot}$ can be defined uniquely.  
 Sufficient for $M_{1} \rightarrow M_{2}$ to preserve algebraicity, is that $\phi(\xi_{1}/\eta_{1})=(\xi_{2}/\eta_{2})$
 a mapping between coefficients, is algebraic. It is necessary, that the scalar product to $d U^{\bot}$
 is mapped projectively on to the scalar product to $d  U^{\diamondsuit}$.
 
Given $d I \in (G_{2})_{C}$, we refer to $d U_{1}=d U_{2}$ as a spiral. Correspondingly, a diagonal in $K(G_{4})$
is represented using $F(x_{j},y_{j})M(x_{j}/y_{j},y_{j}/x_{j})$, j=1,2,3,4, where M corresponds to a symmetric potential
and $(x_{j},y_{j})$ in a convex set. 
In $K(G_{4})$, we get $M=M_{1} \ldots M_{6}$, that is presence of diagonal, implies existence of $M_{j}$ non-algebraic.
The diagonal is characterized by that all of the $M_{j}$ are symmetric, we assume the value of M constant
close to $x_{j}=y_{j}$. For the representation in
$H(K(G_{4}))$, we have that $M_{1}$ algebraic implies $M_{2},M_{3}$ algebraic and conversely. 
Remains the case with 
$x d W + y d W$. We consider $(x^{2} + xy)/(xy + y^{2})=(1 + x/y) / (1 + y/x)$,
that is $M_{4}$ algebraic, if $M_{1}$ algebraic and conversely. The conclusion is that given measures with analytic
coefficients and BV, we have a uniformity space and the polar is given by the diagonals. Note that to conclude
hypoellipticity relative $K(G_{4})$, it is sufficient to exclude spirals in the polar.

Note that for non complete movements, for instance $V \otimes I f(x,y) \rightarrow f(v,0)$, we have that the condition $f(x,y) \leq f(v,0)$
does not preserve relatively compact sets and not complete movements do not give a uniform approximation property,
since the condition on a regular neighborhood is not uniform.

Example: $U - I = I$ implies $Uf=2 I f$, that is invariants to U-I are given by $Uf=2If$. Given f=log g,
and U algebraic, $Ug=g^{2}$. 
Note that $d I \in G(\{ g^{N} \})$ does not imply $d I \in G(\{ g \})$ 
Example: $d I \notin G(\{ \phi \})$, but $d I \in G(\{ \rho \phi \})$
where $\{ \rho \phi \}$ is a nuclear space.

\newtheorem{interpol}[MP]{Interpolation}
\begin{interpol}
Assume $d U = \rho d V$ ($=\rho(v) d v$), for V sub-nuclear, where $\rho$ monotonous on a planar domain $(v,v^{\bot})$. 
Where $\rho$ is independent of $v^{\bot}$, we have $V d V \sim d V^{2}$, that is $(d V) \subset (d V^{2}) \subset (dU)$.
Thus, given $\rho$ linear and independent of $v^{\bot}$, there are intermediate measures.
\end{interpol}

Starting from a very regular boundary and $(d V) \subset (d V^{j}) \subset (d U)$, with $d U = \rho d V$ and 
$d V^{j}=\vartheta^{j-1} d V$,
we can determine dV, so that $\vartheta^{j}$ polynomial.

Example:
Consider two diagonals with respect to $U,V$, $\sigma \sim (V-V^{\bot}) U^{\bot} + V^{\bot}(U^{\bot} - U)$,
that is on den double diagonal $\sigma \sim 0$. Consider $A \rightarrow B \rightarrow C$. Assume $U^{\bot} \sim U_{+}$ and $U \sim U_{-}$ and
in the same manner for V, that is $U_{+}(A)=B_{+}$ and $V_{+}B_{-} =0$,$V_{+}U_{+}=C_{+}$, we then have $\sigma \sim 0$.
Thus, given action for U,V one-sided, there is an approximation property.

\subsubsection{Regular approximations}
For $X_{reg}$ dimension is well-defined. Given $X_{sng}$ removable, we have a continuable movements.
Assume $UW \rightarrow I$ is regular, but $U-I \nrightarrow 0$ regularly, that is $W \rightarrow U^{-1}$.
Spectrum for the approximation, is given by W=I, the character of W is dependent of
type of singularities. Given singularities in a bounded set, then finite iteration of monotonous movements
give a regularization.

The order of G, that regularizes (preserves) singularities, is dependent of the order of singularities.
Assume F a figure for singularities. Example: $V,W$ pseudo local, we then have if $Vf \in C^{\infty}$, that
$WV f \in C^{\infty}$. The condition $WV=VW$ gives, when $V \bot W$, that the order is 2. 

Singularities in H are given by $d Uf=Uf=0$, for instance $Uf \in C \backslash H$. 
Assume U not surjective, but we have existence of W, so that $UW$ surjective. 
Given $R(W^{\bot} U^{\bot})=R(UW)^{\bot}=\{ 0 \}$ (algebraic), we see that $U^{\bot} \bot W^{\bot}$. Further
where $W^{\bot}=(I-W)$ and $U^{\bot}=(I-U)$, we have $U^{\bot} + W^{\bot}=I$ implies $2I - U -W = I$ implies $I = U+W$.
Assume $UW$ a projection operator. Given G has order 2, we must have $W^{\bot}=U$ or $U^{\bot}$,
where $U \neq U^{\bot}$, otherwise we have in the plane 8 possibilities (\cite{Lie91} ch. 4, ch. 17).

\subsubsection{Completion and TP}
Starting from (\cite{Lie91} chapter 5) that two equations A=c and B=$c'$ represent the same set of curves iff
B is only a function of A, $B(x,y)=\Omega(A(x,y))$, we discuss reduction of movements to base movements.

\newtheorem{trivialTP}[MP]{Independence of involution}
\begin{trivialTP}
Consider a regular continuation $d I \rightarrow d V \times d V^{L}$, that is completion to closedness and $d V \times d V^{L} \rightarrow d U$ linear,
this can be given as independent of involution and is not dependent of the type for d V. Assume $d V \times d V^{L}=\mu d I$ and 
$(v,v^{\bot})$ dense in the domain for $1/\mu$ regular, then the completion can be
given as regular.
\end{trivialTP} 
Consider $1/\mu (v,v^{\bot}) d v \simeq d I$. Example: $(G_{2})_{C} \ni d I \rightarrow d I \in (K(G_{4}))_{H}$ completed
to $d V \times d V^{L} \in (K(G_{4}))_{H}$.
Given 
$d V \times d V^{L} \rightarrow d U$ single valued in the plane, then d U can be given as linear in $(d V,d V^{L})$. 
Given an elastic model, $d N \in G$, so that
we have existence of $\rho,\vartheta$, with $\rho d N + \vartheta d N^{\bot}=d I$, then d U can be 
determined as independent of d V, that is we have $d I \rightarrow d U$ regular.  
Assume $d U_{L}=d I_{L}$, where L denotes an entire line,
and that $(d N,d N^{\bot}) \rightarrow (d I,0)$ is rectifiable. More precisely, given
$d \rho / d n^{\bot}=0$ and $d \rho / d n \geq 0$, close to the boundary. Assume $\Sigma=\{ d U =d I \}$, we then have that as long as 
$d V \times d V^{L}$ preserve $\Sigma$, the type for U is not dependent of $d V \times d V^{L}$. 
\newtheorem{tvaspegel}[MP]{Normal model}
\begin{tvaspegel}
Assume $(N,N^{\bot})$ surjective on X, we then have $N \varphi_{1} + N^{\bot} \varphi_{2}=f$, for $f \in X$.
Assume $\phi : \varphi_{1} \rightarrow \varphi_{2}$, that is we have $W \varphi_{1}=f$, where $d W$ is not necessarily 
in the group. Given an approximation property, $N \rightarrow I$ regularly, implies $N^{\bot} \rightarrow 0$
and d W with single valued (monotonous) representation in $(N,N^{\bot})$, we have that $\phi$ does not affect 
the type for d W. Given $N^{\bot} \phi \rightarrow 0$ regularly, we have the approximation $W \varphi_{1} \rightarrow f$
uniformly, without invariant sets. 
\end{tvaspegel}

Given $\psi : d N \rightarrow d N^{\bot}$, we have $\psi^{*}=\phi \rightarrow 0$ close to the boundary.

Assume $d N^{\bot}=\rho d I$, where $\rho$ can be assumed multivalent. Assume $\phi : \varphi_{1} \rightarrow \varphi_{2}$
continuous. Assume $d \phi=\mu d I$, then $\phi$ can be related to a normal system, $d N^{\bot} = (\rho / \mu) d \phi$. 

Consider the extended system and the boundary to the domain $\Omega$, 
given by for instance $(V,V^{L})$.  Assume $d U_{j}=\rho_{j} d U_{1}$,
that is $\rho(u_{1},\ldots,u_{k}) \simeq \rho_{1}(u_{1}) \times \ldots \times \rho_{k}(u_{1})$. 
Note that for $\rho \in C^{\infty}_{c}$, it is sufficient to consider $\rho_{j}$ polynomial. 
Given $\rho \in L^{1}$ and $\rho \rightarrow 0$ in all directions, we have $\rho_{j} \rightarrow 0$ 
$\forall j$. When the domain for $L^{1}$ is bounded, every monotonous movement leaves the domain 
under iteration.

\subsubsection{Decompsable sets}

Consider $d v \rightarrow (d v,d v^{L}) \rightarrow (d u_{1},d u_{2})$. 
Assume p a semi-norm to $(v,v^{L})$ and q a semi-norm to $(v^{\bot}, v^{L,\bot})$.
An extended plane that is contractible to a domain with an approximation property, has an approximation property.
 Presence of a polar set implies that $R(G)$ is not closed. 
 The condition $r \leq pq$ implies that it is sufficient that a semi-norm=0,
 Define $R(\Phi(G))^{\bot} = \{ (w,w^{\bot}) \bot (u,u^{\bot}) \}_{2}$.
$r \leq p q$ defines the polar set, for instance q=0. 

Consider an elastic barrel : $T_{ac}(dU) =\cup_{\mid \rho \mid \leq 1} \rho d U$, with $\rho \in L^{1}$, that is continuation
of the group is defined through $\rho d U$. A necessary criterion for inclusion between weighted $L^{p}-$ spaces 
is for instance $\rho \rightarrow 0$, note that we do not necessarily have existence of a subgroup $(dW)=(T_{ac}(d W)) \subset (d U)$. 
On the other hand, if 
$d U=\sigma d I$, with $\sigma \rightarrow 0$, we can determine d W maximal,
that is a regular subgroup d W, that can be continued as above, to $T_{ac}(d W)=G$.

Given $G_{8}$ generates $R(G)^{\bot}$, we consider $G_{4} \times G_{4} \ni (d V,d V^{L})$. 
Regularity for $(d V,d V^{L})$ in $G_{8}$, corresponds to a two-sided limit.
When we consider $d V - V^{L}=0$,
zero-lines are mapped onto symmetry-sets, $d U - d U^{\bot}=0$ maps zero-lines onto spiral sets.

The condition $(UF-F) + F-I \in C^{\infty}$, implies $(UF-F) \sim (F-I)$ (modulo $C^{\infty}$),
in particular over invariant sets, $F-I \rightarrow 0$ (modulo $C^{\infty}$). Consider 
$(d U)=(d V,d V^{L})$, we then have that invariant sets for d U, corresponds to a domain for UF very regular,
given F very regular. Assume $(v- v^{L})$ is dense in the domain for F and consider $F(v^{L} - v)-I(v^{L} - v) \sim 0$,
that is $(dV,d V^{L})$ is symmetric around I iff $(v- v^{L})$ is symmetric around 0.

Assume $S=\{ f = d f = 0 \}$, that is given $d x \bot d y$, we have $f_{x}=f_{y}=0$.
Note that $f_{u} u_{x}=0$ does not imply $f_{u}=0$, for instance u constant in x.
Assume $\Sigma_{1}=\{ f_{u} \neq 0  \quad f_{x}=0 \}$, in the same manner for $\Sigma_{2}$, $\Sigma_{3}=\{ u f_{u} \neq 0 \quad x f_{x}=0 \}$
for instance $x/u \rightarrow 0$, that is u is reduced with respect to x. In the same manner for $\Sigma_{4},\Sigma_{5}=\{ f_{v} \neq 0 \quad y f_{x}=0 \},\Sigma_{6}$. 
Starting from
$dW = x f_{x} + y f_{y}=0$, with $v f_{v} \neq 0$. $\Sigma_{7},
\Sigma_{8}$ have the corresponding conditions
$x^{2}/v \rightarrow 0$, $y^{2} / v \rightarrow 0, xy /v \rightarrow 0$.
 
 Assume $\Phi(f)=\int f d U$, for instance $\Phi : L^{2} \rightarrow L^{2}$ linear continuous, 
 more generally $\Phi : \mathcal{H} \rightarrow \mathcal{H}$, where $\mathcal{H}$ is not nuclear.
 In particular $\Phi(f) \rightarrow f$, with $d U \notin G_{\mathcal{H}'}$.
 Example: $d I \in G(d u,d u^{\bot})$ but $d I \notin G(d x,d y, d z)$. 
  Presence of non-regular approximations, does not necessarily influence hypoellipticity, that is hypoellipticity does not imply that all approximations are regular. 
 Example: scaling 
 does not necessarily affect P-convexity, however absolute continuity is not preserved (\cite{Collingwood66}, Theorem 3.4), for $d U,d U^{\bot}$ separately.
 Example: if we consider $d U=d I$ iff $d V_{1}(x) \otimes d I(y)=0$, we do not have uniformity and the corresponding
 group is not P-convex, for an iterate of the symbol.

Starting from $\{ f < \lambda \} \subset \subset \Omega$, we consider regular continuations of these
sets, that is we do not have existence of L, an entire line $\subset \{ f=const \}$. Existence of L implies
that convexity (pseudo convexity) is not preserved.

Continuation through continuity allows an entire line in $\{ f=const \}$, that gives unbounded sub-level sets.
Assume $d UW=d I$, where $d U \in G_{2}$. Given $d I \in G_{2}$, we have $d W \in G_{2}$.
Assume $d I \notin (G_{2})_{H}$, we then have $d U \in (G_{2})_{H}$ implies $d U^{-1} \in (G_{8})_{H}$. Assume 
$\Sigma(d U)=\{ (\xi,\eta) \sim (1,1) \}$, gives the coefficients to dU. Given $\Sigma(d U)=\Sigma_{1} \times \Sigma_{2}$
where $\Sigma_{1}=\{ \xi \sim 1 \}$, that is decomposable constant surfaces, we have that invariants to $G_{8}$,
can be derived from invariants to $K(G_{4})$. Note that $(G \times G)_{\mathcal{O}_{M}} \simeq G_{\mathcal{O}_{M}} \times G_{\mathcal{O}_{M}}$,
that is we require slow growth measures.

\subsubsection{Preservation of convexity}

\newtheorem{T/S}[MP]{Two-mirror-model}
\begin{T/S}
A movement, that can be reduced to a normal model, rectifiable at the boundary, can 
be reduced to a single scalar product.
\end{T/S}

Consider $A \rightarrow_{\pi_{1}} B \rightarrow _{\pi_{2}} C$, where TP assumes $\pi_{1}(A) \simeq \pi_{2}^{-1}(C)$.
Example: $d U = \rho d V \rightarrow \rho \mu d I$, where $d I \rightarrow d V$ 
is completion to closedness. For TP, it is sufficient that $1/\rho$ is regular. 

Assume $P_{1} \leq F \leq P_{2}$, where $P_{j}$ are polynomial, then
the sub-level surfaces for F can be given as semi-algebraic. Consider $\Omega=\{ F < \lambda \}$, with $\Omega \rightarrow bd \Omega$
continuous, where the boundary is assumed algebraic. 

Assume $S_{j}(v)=(v,{}^{t} v)=(v,L_{j}(v))$, where $L_{j}(v)$ defines ${}^{t} v$ with respect to the corresponding scalar product. 
Assume $A(v,{}^{t} v)=v + {}^{t} v$. that is $S_{2} A S_{1}=(u_{1},u_{2})$. 
Thus, given the translates are dense in the domain,
$A(v)$ is semi-algebraic, given $S_{j}(v)$ are semi-algebraic. 
Assume $p_{1}(v,v^{\bot}) \leq S^{2}(0) \leq p_{2}(v,v^{\bot})$, where $S^{2}(v,v^{\bot}) \sim S^{2}(0)$
(homotopical conjugation) and $p_{j}$ polynomial. Given isolated singularities and a separation axiom, 
we have that constant surfaces to $p_{j}$ are removable. (Vp=Ip, where I is evaluation in 0). 

Consider the two-mirror-model, when $d {}^{t} V \simeq q_{1} d I$ and $d V \simeq q_{2} dI$,
$d V^{\bot} / d {}^{t} V=p$ and $d V \rightarrow d V^{\bot}$
according to $q_{2} \rightarrow 1/q_{2}$, we then have $d V^{\bot} / d V \sim 1/q_{2}^{2}$, 
given $p q_{1} q_{2} \equiv 1$. We assume $d {}^{t} V \rightarrow d V^{\bot}$ algebraic with an approximation property.

Assume $T(\phi)=f$, where $d U=\rho d I$, where $\rho=1$ describes a line. Assume $d U^{\bot} \leq d W \leq d U$,
where $d U = \rho d I$ and $d U^{\bot}=(1/\rho) dI$ and $1/\rho \leq \sigma \leq \rho$. Cf. Nullstellensats:
$\mid (1/\rho) \rho \mid \leq C$ implies $\mid \rho \mid \leq C'$ for large arguments, given $1/\rho$ reduced,
that is locally 1-1, or equivalently  $\mid d U  d U^{\bot} \mid \leq C$ implies $\mid d U \mid \leq C'$. 
Further, $\{ \sigma \rho < \lambda \} \subset \{ (1/\rho) \rho < \lambda \}$
that is $ \sigma \rho \leq \lambda$ implies $\rho \leq  \lambda$. Conversely, $1/\rho \equiv 1$ (positive measure)
implies $d U = d U^{\bot}$ has positive measure. Assume $U(N,N^{\bot})$ with $d N^{\bot} / d N \rightarrow 0$,
where we assume $(N,N^{\bot})$ has no invariant sets (non-trivial). Thus, $\frac{d}{d N}U(N,N^{\bot}) \simeq d U / d N$
Assume for instance $d U = \rho d N$, with $\rho \geq 1$, that is an one-sided approximation.

\newtheorem{subnuclear}[MP]{Subnuclear movements}
\begin{subnuclear}
Consider $d V^{2} \sim V d V$, given V sub-nuclear, we have $d V=0$ implies $d V^{2}=0$. 
Assume the singularities are given by
$\{ Vf=dVf=0 \}$ and that $d U = \rho d V$, we then have that $\rho=0$ does not necessarily preserve singularities. 
Conversely, given $Uf = WV f$, with W sub-nuclear, we have that regular points are preserved.
\end{subnuclear}

Consider V sub-nuclear  modulo $C^{\infty}$, that is V pseudo local. 
Given $V^{j}=I$ implies $t=0$ (parameter), there are $V^{-1}$ pseudo local, that is over f
the singularities are of finite order.

Given $(U,V) \rightarrow (I,0)$
a nuclear mapping, we do not have a nuclear topology. Sufficient is that I sub-nuclear  (\cite{Schwartz57}).
We assume it is necessary for a nuclear topology,
that we have an approximation property in a convex neighborhood of a regular approximation.  

\subsubsection{A strict condition for TP}

  Consider $f(\mid u \mid) - \mid f(u) \mid \sim 0$. $\mid u_{j} \mid \leq 1$,j=1,2
 corresponding to a cylindrical domain, $\mid f(u) \mid \leq 1$ corresponding to an analytic polyhedron. 
 Let $\Omega=\{ \mid f(u) \mid \leq \lambda \}$ and $\tilde{\Omega}=T_{c}(\Omega)$ (convex closure), corresponds to $(u_{1},u_{2})$
 being dense in $\tilde{\Omega}$.
 
 Consider the condition (\cite{Ramis92}) $e^{B \mid x \mid^{k}} \mid f(1/x) \mid \leq C$ locally. 
 Given $P(1/ \mid x \mid^{k}) \sim 1/Q(\mid x \mid^{k})$,  we have that the condition is $\mid f(x) \mid \leq C Q(\mid x \mid^{k})$,
 note that P polynomial, defines a domain of holomorphy, where s=1/k denotes the size. 
 Given $P(1/\mid x \mid^{k}) \mid f \mid \leq C$, we then have $\mid y \mid / \mid x \mid^{k} \leq C
$ corresponds to $\mid y \mid^{s} \leq C \mid x \mid$ ($ C \rightarrow 0$, $x \rightarrow \infty$). 

\newtheorem{RamisTP}[MP]{TP with strict condition}
\begin{RamisTP}
Assume $R(1/\mid x \mid^{k}) \mid f(x) \mid \leq c$. Assume
$\frac{Q}{P} f(x)$ analytic, with $R(1/\mid x \mid) P(\mid x \mid) \leq c'$, then $\mid f(x) \mid \leq P(\mid x \mid^{k})$ gives a
majorisation condition corresponding to slow growth in $\infty$.  The condition Q reduced
implies $1/Q \rightarrow 0$, that is a strict condition and through Nullstellensats
$\mid \frac{Q}{P} f \mid \leq  c$ implies $R(1/\mid x \mid^{k}) \mid f \mid \leq c'$.
\end{RamisTP}

Concerning $\mathcal{O}_{M}(G_{8}) \simeq \mathcal{O}_{M}(G_{4}) \times \mathcal{O}_{M}(G_{4})$, we
consider $f \in \mathcal{O}_{M}$ according to $\mid f(x) \mid \leq C Q(\mid x \mid^{k})$, for a polynomial Q.
The condition defines a class of measures $\mid Q d U \mid$ BV and immediately regular neighborhoods
of invariant sets. A sufficient condition is a very regular boundary in $(u,v)$.  Assume f corresponding to real type, 
that is $\mid f(v) \mid \leq C e^{B \mid u \mid}$, with B close
to 0, then the invariant sets can be determined in $K(G_{4})$.

\subsection{Projective transformations}
\subsubsection{Topological degree}
Assume a normal model $(N,{}^{t} N) \rightarrow (N,N^{\bot})$, that is $N \rightarrow I$ iff $N^{\bot} \rightarrow 0$.
When N is algebraic, ${}^{t} N$ is algebraic. Assume ${}^{t} N \rightarrow N^{\bot}$ completion to
an approximation property. Sufficient for an algebraic continuation is strict pseudo convexity.
Algebraicity for d N is dependent of topology. 
Example: $d N(\tilde{\phi})=d \tilde{N}(\phi)$, where $C^{\infty} \ni \phi \rightarrow \tilde{\phi} \in L^{1}$ with 
$d \tilde{N}$ algebraic. 
Starting from $N(U,U^{\bot})$, where U one-parameter, a sufficient condition for a single valued orientation close to the boundary
is $d U^{\bot} /  dN \rightarrow 0$. For hypoellipticity, it is necessary that $N(U,U^{\bot})$ is corresponding to real type,
that is determined by d N / d U.

Projective transformations are the analogue to transformations that map closed forms onto closed forms.
For regular transformations, we assume that a closed contour is mapped on a single Jordan curve in the plane, 
or a cycle in $S^{n}$.
Given a strict condition, there is to a projective mapping, a corresponding bijective mapping.
Example: $d U = p d V$, where p is reduced, we have $\int_{L} d U =0$ implies $\int_{L} d  V=0$ or $L = \{ 0 \}$.
A strict homomorphism, corresponds to a local bijection $dV \rightarrow dU=\rho d V$. 
Sufficient for isolated singularities, is $d U=0$ implies U=I. Note that d U reduced does not imply
U reduced. Since every cycle in $S^{n}$ can be mapped on to a full line in $E^{n+1}$, 
it is to determine topological degree (\cite{Leray34}), sufficient ot consider zero- lines for measures. We consider homotopy,
represented as $d U = \rho_{j} d N$
with $\rho_{1} \sim \rho_{2}$, that is given an approximation property, the degree is preserved.

\newtheorem{grad}[MP]{Order for projective transformations}
\begin{grad}
Assume cycles $\rightarrow$ zero-lines, we then have that the topological degree gives the order for projective transformations.
\end{grad}

Consider a zero-line L to $d U$ BV, that can be mapped on $\Gamma$, a simple Jordan curve. 
Assume further $L'$ a zero line to $d U^{\diamondsuit}$ and to $d U^{\bot}=\sigma d U$, then d U=0 on $(L,L')$, 
that is all zero-lines between $L,L'$,
that can further be mapped onto cycles (\cite{Collingwood66} Lindelöf's theorem). Assume $d U = \rho d I$ and $d U^{\bot}=\rho^{\bot}d I$, 
with $\rho - \rho^{\bot} \in H (L,L')$, we then have given $L \rightarrow L'$ bijective, that $\rho-\rho^{\bot} =0$ on $L,L'$, implies
$d U = d U^{\bot}$ on the sector between $L,L'$. Given instead $d U = d U^{\bot}$, that is $H \ni \rho \rightarrow \rho^{\bot} \in \mathcal{H}$ and $\rho=\rho^{\bot}$ on $L,L'$ and that 
we have existence of a radius $\gamma$ on a convex set between $L,L'$ (planar), so that $d U \neq d U^{\bot}$ on $\gamma$, 
then, given H dense in $\mathcal{H}$, there are no spirals (diagonal elements) in the inner of the domain between $L,L'$.
Consider $d I \rightarrow d N \rightarrow d U$, where d N is a boundary operator.
We assume $d N^{\bot} \rightarrow 0$ regular, thus we must have for the polar set $d U^{\bot}=\sigma^{\bot} d N^{\bot}=0$, 
that $\sigma^{\bot} \equiv 0$.
Consider $Y(f) = \{ u^{\bot} \quad f(u,u^{\bot})=0 \}$ a manifold for the polar.
Given that the measure for $Y(f)$ is non-zero, then when $\rho \rightarrow 1$, we have multivalentness for the limit (\cite{Collingwood66}). 
Let $\phi_{j} : d U_{j} \rightarrow \rho_{j} d U_{1}$, where 
$\rho_{j} \neq 0$ define continuous mappings. The order for Y, can be given as the number of j, such
that $\delta \sigma^{\bot} / \delta u_{j}=0$.

Consider the $T_{ac}$-closure, $T_{ac}(d V)(\mathcal{H})=\{ \rho d V \quad \rho \in \mathcal{H} \}$, that is
$d V \rightarrow \rho d V$ is projective. Given $d V= \mu d I$, $\mu \rho= \mu + \rho$
iff $1/\mu + 1/\rho=1$, that is the condition implies a dual multiplier. 
Given $\rho(v,v^{\bot})$, where $\rho$ single valued, then $\rho=\alpha \rho_{1}(v) + \beta \rho_{2}(v^{\bot})$
and $\rho= \alpha \rho_{1}(x)\mu(x) + \beta \rho_{2}(x) \mu^{\bot}(x)$.  
Example: Assume conjugation according to $d U = d I$ iff $d U^{\bot}=0$, that is $d U^{\bot} \in G$ iff $d U - d I \in G$
but $d I \notin G$. A necessary condition for $d I \in G$, is a regular neighborhood of $ d U=d I$. 

Assume $\Phi(f)=\int f d U$, assume $f=T_{x}(\phi)$, why given that the topology for $\Phi$ is nuclear,
we have $<T_{x},d U>(\phi)=\Phi(f)$, where we assume the topology invariant for U. In a non-nuclear
topology, we do not necessarily have $\Phi(f) - f \rightarrow 0$ implies $(T_{x} - I) \rightarrow 0$.
Example: assume $\mathcal{H}_{1} \subset \mathcal{H} \subset \mathcal{H}_{2}=\mathcal{D}'$. 
That is assume $\Phi(f) - f \rightarrow 0$ in $\mathcal{H}$ and 
$(T_{x} - I) \rightarrow 0$ in $\mathcal{H}_{1}$, that is we assume $\mathcal{H}_{1}$ nuclear.
Further, $\Phi(f)-f \rightarrow 0$ implies $T_{x} - I \rightarrow 0$ in $\mathcal{D}'$, 
but not in $\mathcal{H}$.

Assume existence of $\Phi$ (ähnlich), with $dU=\Phi(d I) \in G$, but $\Phi^{-1}(d U) = d I \notin G$.
Consider $K(G)$, that is $(d V,d V^{L}) \in G \times G$. Starting from $G_{8}$,
consider $G_{8}'=\{ d  V \quad d V + d V^{L} \in G_{2} \}$. Consider $G_{8}''=\{ d V \quad \Phi(d V,d V^{L}) \in G_{2} \times G_{2} \}$
and consider $G_{8}^{ac}=\{ d V \quad \Phi(d V,d V^{L}) \in G_{2}^{ac} \times G_{2}^{ac} \}$, where
we assume $d U \in G_{2}^{ac}$ with $d U=0$ implies $U=I$ (invariants). Further, $H(K(G_{4})) \subseteq C(G_{2})$,
where $d I \notin (G_{2})_{H}$, but $d I \in (G_{2})_{C}$. Given $\Phi \rightarrow d I$ regularly, we have an approximation property.
Assume $\{ \Phi < \lambda \}$ is semi-algebraic (and BV). We then have $d I \notin G_{2}$ (in H), 
but $\Phi(d V,d V^{L}) \in G_{2} \times G_{2}$,
that is we consider $G_{8}''$ above, that is not necessarily semi-algebraic in $G_{8}$. 

Consider $f(u,v)$ with isolated singularities, that is $\{ f(u,v)=d f =0 \}$ are isolated. 
Obviously $\{ (u,v) \quad f=d f=0 \} \subset \{ (x,y) \quad f=d f=0 \}$. $f(u,v) \in L^{1}$ 
is related to $f_{0}$ very regular in $L^{1}$ and when we assume $f \rightarrow f_{0}$,
we have $d I \in G_{2}$.

Consider $d V \rightarrow d V \times d V \rightarrow d U \times d U^{\bot}$, in particular where d U
is determined by d V only, ex $d U \rightarrow d V \rightarrow d I$ simultaneously $d U^{\bot} / d V \rightarrow 0$,
for instance d U of type 0, real dominant. Note that $d V \times d V \rightarrow d U \times d U^{\bot}$
can be given as projective, given $d U^{\bot} \bot d U$. Assume that the support for f
includes $\{ d U = d I \}$ and given existence of $d V=0$ over invariant sets to dU (projective conjugation),
we have that a neighborhood of invariant sets can be given as a regular analytic domain. 
Consider $d V(f)=\xi f_{x} + \eta f_{y}$, we then have  $d V(f)=d V(u) f_{u} + d V(u^{\bot}) f_{u^{\bot}}$.

\subsubsection{Partial conditions}

\newtheorem{phe}[MP]{Partial regularization}
\begin{phe}
Assume $f(u_{1},\ldots,u_{k})$ have regular approximations in some directions. Given a condition that 
invariant sets in irregular directions are subsets of  
a domain of holomorphy and given $f$ partially hypoelliptic, $f^{N}$ has regular approximations in all 
directions.
\end{phe}

Assume $d U + d W=d I$ in nbhd $\Sigma(d U)=\{ d U = d I \}$, where dU in $G$.
Thus, the polar set is subset of $d W=0$. Note that $WI=IW$ does not imply $d W=0$. 
Linear phase for W acting on algebraic elements, is sufficient for an algebraic polar set.

Consider $H(G_{2}) \subset H(K(G_{4}))$, with the condition that every domain $\Omega$, has a movement in $H(K(G_{4}))$,
with $\Omega$ as a domain for analyticity. Example: $f \in C(G_{2})$
with $d U_{3}(f)=0$, where $d U_{3} \in K(G_{4})$. Then given $f(u_{1},u_{2},u_{3})$ partially hypoelliptic, 
there are N, with $f^{N} \in H(\Omega)$.

Consider $U \rightarrow U^{\diamondsuit} \rightarrow U^{\bot}$, where the last mapping is trivial.
(\cite{AhlforsSario60}), given $d U=0$ iff $d U^{\diamondsuit} =0$, we have an analytic representation.
Pure mappings preserve analyticity in the plane, that is symmetry with respect to $(x,y) \rightarrow (y,x)$, 
further $(x,y) \rightarrow (x,-y)$. Given $(d U)^{\bot} \simeq d U^{\bot}$, with the conditions above, 
we have an analytic orthogonal. Assume partial P-convexity, with symmetric domain otherwise, 
we then have $f^{N}$ has a corresponding P-convex G, in particular an approximation property 
in all variables separately. 

The condition on invariant sets in a domain of holomorphy, implies $\{ U-I \} \subset \{ d W=0 \}$, where $d W$ analytic,
that is U has an analytic projective conjugate, that is d W with d U + d W=d I, that can be continued to 
$U^{\bot}$. Consider $\pi : (x,y,z) \rightarrow (u,v)$,
so that $\pi^{-1}(0) \sim (x,y,z)$. Example: $(d U, d U^{\bot})$ can be reduced to a normal system $(d N,d N^{\bot})$, 
so that  $d N(U,U^{\bot})$ is only dependent of d U close to the boundary. Alternatively, consider ($\rho d U_{1},\vartheta d U_{2})$
with $\rho d U_{1} \rightarrow d I$ and $\vartheta d U_{2} \rightarrow 0$, as regular approximations.
We can assume, in particular $(d U_{1},d U_{2}) \rightarrow (d I,0)$ in a normal model. We have the case 
$d I \notin G$, when $1/\rho$ is not regular. 

Example: assume $f(v_{1},v_{2}) \rightarrow 0$, when $v_{1} \rightarrow 0$ is regular, that is $v_{1}=0$ 
has a regular neighborhood. In the converse case, assume that
$v_{2}$ can be continuously deformed to a non-regular direction, for instance $f \nrightarrow 0$, 
$v_{2} \rightarrow 0$. Sufficient for a
regular neighborhood, is that f is algebraic in $v_{2}$.  

Example: assume u,v one-parameter, with $f(u/v) \rightarrow 0$ regularly, but not $f(u) \rightarrow 0$ regularly.
If $f(u,v) \rightarrow 0$ regularly, when $u \in L$ and $u \rightarrow 0$, that is L gives a regular neighborhood of 0,
at the same time, $f(u,v) \rightarrow g(u)$, when $v \rightarrow 0$, that is $\lim_{u \rightarrow 0} \neq \lim_{v \rightarrow 0}$.
Then, we do not have regular continuations according to a normal model. However, there may be $d W \in G_{8}$,
with $f(u,v,w) \rightarrow 0$ regularly.
Assume the boundary is given by $U=I$ and that $R(d U)^{\bot}=R(d V)$ and that the polar is defined
by $Vf=0$. A very regular boundary over $K(G)$, can be represented by $f^{(k)}(u,v)=const$
implies $u=v=0$.
For a normal model, we require that 
$\lim_{u \rightarrow 0} \frac{\delta f}{\delta v}(u,v)=0$ regularly.
Note that $f^{N}$ is corresponding to real type, does not imply that f is corresponding to real type.

\subsubsection{Reduction to normal model}

\newtheorem{invers}[MP]{The group for the inverse mapping}
\begin{invers}
Assume $F(\gamma_{T})(\zeta)=F(\gamma)(\zeta_{T})$, then $U \gamma_{T}(\zeta)$ is
linear in $\gamma_{T}$, but only corresponds to a transformation in $\zeta_{T}$.
\end{invers}

Assume $\zeta_{T}$ through a transformation in $C(G_{2})$, simultaneously analytic in $K(G_{4})$, 
that is $G_{2} \gamma (\zeta)=\gamma(K(G_{4}))$.
Note that, for $\gamma(\zeta)=(\gamma_{1}(\zeta_{1}),\gamma_{2}(\zeta_{2}))$,  $\xi f_{\zeta_{1}} + \eta f_{\zeta_{2}}$
$=\xi \gamma_{1,\zeta_{1}} f_{\gamma_{1}} + \eta \gamma_{2, \zeta_{2}} f_{\gamma_{2}}$. Thus,
given a movement in $\zeta$, we can derive a movement in $\gamma_{1},\gamma_{2}$, according to
$\zeta \rightarrow \gamma \rightarrow F(\gamma)(\zeta)=f(\zeta)$.
Assume $F(U \phi_{T})(\zeta)=F(\phi)(\zeta_{T})$, then given d U is analytic, $\zeta_{T}$ is defined
through a continuous transformation. Example: $d U_{1}(f)(\zeta)=\xi(\zeta) f_{x} + \eta(\zeta) f_{y}$.
Given $\zeta \rightarrow (x,y)$ analytic, through Cauchy-Riemann and if we consider $\zeta_{2}$ as parameter,
$\xi=(y+x)_{\zeta_{2}}$ and $\eta=(y-x)_{\zeta_{2}}$. Thus $\xi + \eta$ is only dependent of y and
$\xi - \eta$ is only dependent of x.

Assume $U$ sub-nuclear , that is $d U =\rho d V$, where $d V \in G_{r}$ and $d U \in G$. Given 
$\rho \rightarrow 1$ regular and $\mathcal{H}(d V)$ nuclear, we have $\mathcal{H}(d U)$ nuclear. 
A special case of TP is $U \rightarrow U_{r} \rightarrow \zeta$,
where $\mathcal{H}(d U)$ nuclear. Thus given $d U=0$, we have existence of $d U_{r}=0$, 
that is $\rho \neq 0$ and $1/\rho \in \mathcal{H}$.

Assume $d U(x,y)$ analytic over supp f with $f \in (I)$. We then have
$f(u,u^{\bot})$ analytic implies $f(x,y,z)$ analytic. Starting with $F(\gamma_{1} \rightarrow \gamma_{2})(\zeta)$
analytic, it is necessary that $\gamma_{2} \rightarrow \gamma_{1} \rightarrow F(\gamma_{1})$ continuous.
We consider $\gamma_{1} \rightarrow \gamma_{2}$ as topological completion.

\newtheorem{irreducibler}[MP]{Reducibles}
\begin{irreducibler}
Assume $d N = \rho d U$ a reduction, with $(d N,d N^{\bot}) \rightarrow (d I,0)$ regular and $d N=d I$ rectifiable. 
Then $\rho(u,u^{\bot}) \rightarrow 1$ as an absolute continuous function and the reduction preserves zero-sets (\cite{Riesz20}).
\end{irreducibler}

More precisely $\rho : \{ d U=d I \} \rightarrow \{ d N = d I \}$ is absolute continuous and maps sets of measure
zero on sets of measure zero. We assume $\rho$ analytic, where $N^{\bot} \neq 0$. 
In particular the limit of $\rho$ is not dependent of path for the approximation.
Given $\rho$ monotonous, that is single valued, we have close to the boundary, $\rho(u,u^{\bot}) \sim c_{1} u + c_{2}$, 
that is we do not have a spiral axes.
Example: $d U^{\bot}=\sigma d N^{\bot}$, with $\sigma^{k}$ polynomial,
alternatively $d U^{\bot}=\mu d U$, with $\mu$ polynomial. 
We can obviously determine a domain for absolute continuity in $\{ (u,u^{\bot}) \}$, where $\delta \rho / \delta u^{\bot} \rightarrow 0$
and we have two-sided limits, when $\pm u^{\bot} \rightarrow 0$ regularly, cf sharp fronts (\cite{Garding77}).

Example: $d U = \rho d V = \rho \sigma d I$.
Given $\rho$ bounded and absolute continuous, with $Vf=If$ iff $f(v)=f(0)$ $(=f(x,y,z)$), we have where $d \rho=0$ that $d U=d I$ implies $d V=d I$, 
that is we have reduction of dU to for instance a normal model.
Example: $d U^{\bot}=d V^{\bot} W^{\bot}$, where $d V^{\bot} \rightarrow 0$ regularly. Assume $d W=p d I$
and $d W^{\bot}=1/p d I$. Then we must have $p d U^{\bot} \rightarrow 0$ regularly.  

Example: Assume improperly $VIf=f(v-0)$ and $IVf=f(v+0)$,
given $VI=IV$ over the support to f and given one-sidedness for R(V), we do not have traces. 
We conclude that one-sidedness for R(V) preserve convexity in the plane.
An irreducible continuation $d U = \rho d V$, is unique. When $\rho$ is linear in V, 
the continuation can be given through iteration. Given $d U=d I$ implies d V = dI, 
the continuation preserves type.
Assume for real numbers $\alpha,\beta$, $f(\alpha u + \beta u^{\bot}) \leq \alpha f(u) + \beta f(u^{\bot})$, 
where $\alpha + \beta=1$, that is f is convex (absolute continuous).
Given $f(u) \rightarrow f(0)$, we have $f(u^{\bot}) \rightarrow 0$. 

Assume $f(u)=f(0)$ implies $u \in L$, for instance f 1-homogeneous over invariant sets. 
Obviously, we have existence of $d V$ conjugated,
so that f (-1) homogeneous with respect to v. We assume f independent of v on L.
For irreducibels, we have $f(u,v)=f(0)$ implies $U=I$ or $V=I$, for instance a projective model,
invariant on L, where $L \sim 0$. 
\newtheorem{volume}[MP]{Volume preserving conjugation}
\begin{volume}
It is not clear that the reduction model above, is independent of 
scaling $f \rightarrow \tilde{f}$. However, given $d U \rightarrow d U^{\bot}$ conjugated according to 
$t \rightarrow 1/t$ (scaling parameter), we have that $\tilde{f}(u,u^{\bot})$ has the same volume as $f(u,u^{\bot})$. 
P-convexity for the group
implies further, that $U^{\bot}f=0$ implies t finite. 
\end{volume}

Consider $f(x,y)M(y/x,x/y) \simeq f(x,y/x) f(x/y,y)$, with M=1, we then have $f(x,y)=0$ implies
$f(x,y/x)=0$ or $f(y,x/y)=0$ (or both). Irreducibility excludes the case both.
The condition M=0 for irreducibels implies $U=U_{1}$ or $U=U_{2}$,
that is a decomposition into irreducibels excludes spirals.
\subsubsection{The polar group}

Example: consider $E$, the localization of f, 
$C=\{ (u,u^{\bot}) \quad E(f)(u,u^{\bot})=0 \}$, where $E-I=0$ modulo regularizing action, 
C gives the polar set to hypoellipticity.
Assume $(u,u^{\bot})$ is dense in the domain for f, E has type 0 over $u=u^{\bot}$ and does not correspond
to a regularizing action.
Assume $E_{k}$ corresponds to the localization of $f^{k}$, where f is assumed partially hypoelliptic, 
that is $f^{k}$ locally hypoelliptic, for large k.
We then have ker $E_{k} \downarrow \{ 0 \}$ (modulo regularizing action), when $k \uparrow \infty$.

Consider $C(G_{2}) \rightarrow H(K(G_{4}))$, where we assume $H(K(G_{4}))$
with an approximation property, that is uniform convergence on convex, compact equilibrated sets.
Assume further a strict condition, that $(v,v^{L}) \rightarrow v + v^{L}$ continuous and surjective.
Then we have that $C(G_{2})$ has an approximation property.

$\Omega(x,y,y',y'')= y'' - \omega(x,y,y')=0$ represents at most 8 infinitesimal transformations in the plane 
(\cite{Lie91}). Maximal for $\omega \equiv 0$. $y''=0$ has as integral curves $\infty^{2}$ straight lines in the plane.
$\Omega=0$ can be deformed to $y''=0$. 
In the plane there is a domain $\Sigma$, such that $p,q \in \Sigma$ implies existence of a unique 
integral curve through $p,q$.  

Example: projective transformations map lines in the plane on lines in the plane. Through Theorem 35 (\cite{Lie91}, cf. ch. 17 prgf. 6)
$-\xi_{yy} \rho^{3} + (\eta_{yy} - 2 \xi_{xy}) \rho^{2} + (2 \eta_{xy} - \xi_{xx}) \rho + \eta_{xx}=0$,
where $\xi,\eta$ is only dependent of $x,y$. Since $\xi_{yy}=0$ implies $\xi=X(x) y + X_{0}(x)$ and in 
the same manner for $\eta_{xx}=0$.
Simultaneously, $\eta_{y} - 2 \xi_{x}=X_{1}(x)$, $\xi_{x} - 2 \eta_{y}=Y_{1}(y)$. Thus $3 \xi_{x}=-2 X_{1}(x) - Y_{1}(y)$
and $3 \eta_{y}=-2 Y_{1}(y) - X_{1}(x)$. Further, $3 X' y + 3 X_{0}' = - 2 X_{1}(x) - Y_{1}(y)$,
$3 Y' x + 3 Y_{0}'=-2 Y_{1}(y) - X_{1}(x)$. The conclusion is, $Uf =(a + cx + dy + h x^{2} + k xy) \frac{\delta f}{\delta x}
+ (b + c x + g y + h xy + k y^{2}) \frac{\delta f}{\delta y}$. Thus, there are 8 independent transformations
that are generated by $f_{x},f_{y}, x f_{x}, y f_{x}, x f_{y}, y f_{y}$,$x^{2} f_{x} + xy f_{y}$ and
$xy f_{x} + y^{2} f_{y}$.

 Note $Wf=g$ with $Wf = \int f_{1} d W$, for some $f_{1}$, where g=$Ig =\int g d I$, that is solvability
 does not imply an approximation property for the solution.
 The condition that $d U=0$ over $\gamma_{1},\gamma_{2}$, implies that the mapping
$\gamma_{1}(u) \rightarrow \gamma_{2}(u)$ is projective. That is dU defines a geometric ideal, $u \in N(I)$.
Assume $L_{1},L_{2}$ are invariant lines to $d U_{1},dU_{2}$, further $L_{1} \rightarrow L_{2}$
projective, we then have if $d U_{1}=0$ over $L_{2}$, $\gamma_{1}(0)=0$ implies $\gamma_{2}(0)=0$.

\subsection{TP and P-convexity}
Given strict pseudo convexity, we have $H(G_{2}) \subset L^{1}$, with an approximation property, 
in the converse case,
we have that $H(K(G_{4}))$ approximates $L^{1}$ uniformly.

We consider the transmission property over $T_{ac}(dI)=\cup_{\mid \rho \mid \leq 1} \rho d I$, 
where  $\rho \rightarrow 1/\rho$ have the same constant surfaces, that are assumed trivial. 
Example: $\rho = p / q$, where either
$p/q \rightarrow 0$ in $\infty$ or $q/p \rightarrow 0$ in $\infty$.

Assume for every $d U(x.y.z)$, there are $\Omega$, so that $d U$ are analytic over $\Omega$. 
Consider $d U = \rho d I$, given $\rho \in L^{1}$, given a strict condition, there are $\rho \in H$ 
for which $d U$ is analytic. Assume $\phi \rightarrow \phi_{1}$
translation of the support to $\Omega$, so that $d U$ is analytic over the support to 
$\phi_{1}$, then modulo translation, we can assume dU analytic. 

Consider $U \rightarrow {}^{t} U \rightarrow U^{\bot}$, under a strict condition, 
with $d U^{\bot}=0$ (polar set). The TP model (special case) gives
$F(\gamma_{1} \rightarrow \gamma_{2})(\zeta) \rightarrow \zeta$, where $\gamma_{2}$ analytic, 
can now be continued to $\mid \zeta \mid^{\sigma} \leq C \mid \gamma_{3} \mid \leq \mid F(\gamma_{1}) \mid$, 
that is $\gamma_{3}$ is locally reduced. Sufficient is $\mid \zeta \mid^{\sigma} \leq C\mid \gamma_{3} \mid^{N} \leq \mid F(\gamma_{1}) \mid$, 
that is a sufficient condition is $\gamma_{3}$ partially hypoelliptic. In this case the polar set can contain spiral sets, 
that is we do not have a strict condition. Thus, given TP, we can prove that P-convexity
implies hypoellipticity.

Topological singularities: Consider $\gamma_{1} \rightarrow \gamma_{2}$ according to $\mathcal{H} \rightarrow H$
or $C^{\infty}$, that is we have completion to $C^{\infty}$
through a subgroup $G_{reg}$. Note that $T(\phi) \in C^{0}$ is necessary
for a representation with a non-decreasing Stieltjes measure.

Assume $h(f/g) \rightarrow 0$ defines convergence in norm, that is an approximation property with respect to norm.
We assume existence of W, so that $W \mid f \mid \simeq f$ (rotation), relative topology. Thus, convergence in norm with respect to U, implies mean convergence with respect to U,W. 
W is assumed pseudo local with $W \rightarrow I$ in $C^{\infty}$, thus (topologically) algebraic. 

\newtheorem{lyftprincip}[MP]{Condition for mean convergence}
\begin{lyftprincip}
 Assume $d U = \rho d V = \rho \mu d I$, where $1/\rho \notin \mathcal{H}'$. Assume $h(f/g)=0$ implies $\mid \mid f \mid - \mid g \mid \mid \rightarrow 0$
that is ker h defines a symmetric ideal (complex conjugation). However, $\mid g \mid$ polynomial, 
does not imply that f is approaching a polynomial, but given isolated singularities for h, 
we can assume $f \rightarrow g$ in the mean.
\end{lyftprincip}

Existence of d W, so that $W(\mid f \mid - \mid g \mid) \sim \mid f- g \mid$, is motivated by 
the convergence being preserved under $S : (x,y) \rightarrow (x,-y)$ and ${}^{t} Id : (x,y) \rightarrow (y,x)$. 

Assume $d UV = \rho \mu d I$, where $d U \in C(G_{2})$ and $d V \in H(K(G_{4}))$. Sufficient for $f \rightarrow \mu f$
to preserve analyticity, is that $\mu$ is analytic. Given $1/\mu \sim \rho$, we require that $f \rightarrow \rho f$ preserve
continuity. Example: $d \mu \neq 0$ (conformal) gives $\rho \in C$. We have isolated singularities for h, only if
the limit involves both $U_{1}$ and $U_{2}$, that is completeness for the limits. 

\subsubsection{Representations}

Assume $d U \rightarrow d V$ contractible, in particular a restriction
without changing type, that is without affecting invariant sets. Note that two scalar
products are involved.  
Consider $(d U,d U^{\bot}) \rightarrow (d U, d U^{\bot} / d U) \rightarrow (d I , 0)$, defines a
weak orthogonal in $\infty$, corresponding to a strict condition for convergence. Alternatively,
determine $\Gamma$ with $d N \neq 0$ on $\Gamma$ and $d N \bot d U$, that is dU is a continuation of d N, 
analogously with P-convexity.
Simultaneously, $d N \rightarrow d I$ on a component (contractible).
Note, using Schwartz kernel theorem, $<<T(\phi),d U>,d U^{\bot}>=<T(\phi), <d U,d U^{\bot}>>$. 

Assume $C(G_{2}) \backslash H(G_{2}) \simeq H(K(G_{4}))$, where H approximates $L^{1}$, given 
strict pseudo convexity. The condition $d U \bot d U^{\bot}$, is a weak condition with respect to $C(G_{2})$, 
that is it is defined by Stieltjes measures. The completion to $L^{1}$ defines $d U^{\bot}$ (\cite{AhlforsSario60}).
Consider $<T_{x.y}, V^{L} \phi_{1} \otimes V \phi_{2}>=<T_{x,y}, d V^{L} \otimes d V>(\phi_{1} \otimes \phi_{2})$.
Assume $T_{x,y}$ very regular and consider $T(\phi_{1} \otimes \phi_{2})=0$, that is movements $d V^{L} \otimes d V$,
that preserves a zero space. Thus, $d V^{L} \otimes d V$ describes the orthogonal to T over $\phi_{1} \otimes \phi_{2}$.

Assume $d U d U^{\bot}=\mu d N d N^{\bot}$, where the boundary is defined according to P-convexity. 
Given $d U d U^{\bot}$ gives the volume for $\Omega$, $d N d N^{\bot}$, when $N^{\bot} \rightarrow 0$, 
gives the measure for $bd \Omega$. A necessary condition for TP is that
$\Omega \rightarrow bd \Omega$ continuous, that is we assume $d U d U^{\bot}$ invertible, that is $\mu \neq 0$.
Example: assume $d U = \rho \sigma d I$, with $d U \rightarrow d N$ regular, that is consider $\rho \sigma$ regular. 
Given $d U =\rho d N$ and $d U^{\bot}=d N/\rho$, we have $d N(U,U^{\bot})=(1/\rho) d U + \rho dU^{\bot}$,
where $\sigma / \rho \geq 0$ close to the boundary.
   
 Consider $U (N,N^{\bot}) \rightarrow N(U,U^{\bot})$. When $N \rightarrow I$, $U,U^{\bot}$ describes the polar set. According to Gauss-Greeen $- \int_{S} - Yd U- X d U^{\bot}=\int \int_{A} (X_{U} + Y_{U^{\bot}}) d U d U^{\bot}$.
 Assume for instance Y,X are given by $d N / d U$, $d N / d U^{\bot}$. Given N symmetric in $U,U^{\bot}$, 
 we have $\int_{S} d N=0$, when $N \rightarrow I$ regularly, why A must be trivial in a normal model.

Example: $d U \rightarrow d I$ with $d U / d N^{\bot} \rightarrow 0$,
simultaneously $d U^{\bot} \rightarrow 0$. An algebraic continuation preserves reducedness, 
in particular $d N \rightarrow d U$ is reversible, when $d N^{\bot} / d N \rightarrow 0$.

\subsubsection{A transmissions group}

Starting from the condition, that $U \rightarrow {}^{t} U$ preserves type, we can assume a Dirichlet problem,
that is $\eta_{x} - \xi_{y}=0$ ($\xi_{x} + \eta_{y}=0$), for $\mathcal{O}_{M}(G_{8}) \simeq \mathcal{O}_{M}(G_{4}) \times \mathcal{O}_{M}(G_{4})$ (\cite{Schwartz57}). 
Consider  $A \simeq d V \times d V^{L} \rightarrow d U \times d U^{\bot}=C$ BV, as a continuous continuation. 
Assume $B=d N \times d N^{\bot}$ defines a rectifiable boundary, with $N^{\bot} / N \rightarrow 0$ and $d N \rightarrow d I$. 
With these conditions, we assume TP is given by $U \times U^{\bot} \rightarrow N \times N^{\bot}$ continuous (absolute continuous). 
With the condition that
$d N \times d N^{\bot}$ is P-convex, we have that hypoellipticity is preserved. Sufficient for preservation of $T_{ac}$-convexity,
is that $dU \rightarrow 1/\rho d U$ is absolute continuous.

\newtheorem{TPsubgroup}[MP]{A TP subgroup}
\begin{TPsubgroup}
Assume Y : $\gamma_{1} \rightarrow \gamma_{2}$ projective on test-functions and for instance $F(\gamma_{1})$ continuous, but 
$F(\gamma_{2})$ analytic. Example: $A \rightarrow B \rightarrow_{Y} C$
and $A \rightarrow_{U} C$, that is $Y^{-1}U : A \rightarrow B$ and $U^{-1} Y : B \rightarrow A$.
\end{TPsubgroup}

Note that given strict pseudo convexity, we can define $L^{1}(bd X) \rightarrow H(bd X)$, 
that is a monotropic operator, further the scheme above can be given as symplectic.

Assume F defines singularities and $WF \cap F = \{ 0 \}$, that is dW is in the regularity group, 
then G can be defined as irreducibels relative W, that is $U \neq W^{-1}$ over F. 
Given $d W^{2}$ has isolated singularities, when $Wf \in C^{\infty}$ for $f \in C^{\infty}$, 
we have that W can be represented as algebraic, analogous with an approximation property.

Consider $d V \subset d U \subset d V$, with $d U = \rho d V$ and $\mid \rho \mid \leq 1$,
that is $d U \in T_{ac}(d V)$. We assume $\rho \rightarrow 1$ at the boundary,
further $d V=\sigma d U$, with $\sigma \rightarrow 1$ at the boundary. Example: $d U=d V$ outside compact, 
and otherwise one-parameter $p d U=q d V$, for polynomials p,q, that is we have $p/q \in H$,
except for finitely many terms. Example: assume $d U^{\bot}=\rho^{\bot} d V^{\bot}=\rho^{\bot} \mu^{\bot} dI$ and that $\rho^{\bot} \mu^{\bot} \rightarrow 0$, as $t \rightarrow \infty$,
simultaneously $\rho \rightarrow 1$, as $t \rightarrow 0$.

\subsubsection{Functional transformations}
P-convexity implies that $U^{\bot}f \neq 0$ in the infinity. Further, assume existence of d N with $Nf \neq 0$
and $d U \bot {}^{t}p d N$ on the boundary. $d I \in G(R(N))$ (range) implies $d UN=\rho d N$, with $1/\rho$ regular. Given
$U(N,N^{\bot})$ single valued, we have $d U/d N=1$ at the boundary, that is $d U \rightarrow dN \rightarrow d I$ regular. 
Note that,
$U^{\bot}\phi=0$ implies $u^{\bot}$ finite, that is $u^{\bot} \nrightarrow 0$, but $u^{\bot}/u \rightarrow 0$.

Starting from $H(K(G_{4})) \subseteq C(G_{2})$, we have in H, that $f \equiv const$ on a boundary of positive measure,
implies $f \equiv const$ on a convex neighborhood. Through the homotopy mapping we can prove, 
given f has a zero in the neighborhood, that $f \equiv 0$,
for instance $f(K(G_{4})) \sim_{0} g(G_{2})$, gives that g reduced to a normal model is $\equiv 0$, and 
can be continued to a convex neighborhood. 

\newtheorem{normal}[MP]{Normal model}
\begin{normal}
Assume $\Phi(u)=u-\mathcal{F}(u)$, for instance $\mathcal{F}(u)=u^{\bot}$ (conjugation), starting from 
a strict condition,
for instance $u \rightarrow 0$ iff $U^{\bot} \rightarrow 0$, we have that $\Phi$ bijective,
corresponds to density for $(u,\Phi)$.
\end{normal}

Assume $\Phi$ defined by $d \sigma$, for instance Dirichlet measure, L a zero-line to $d \sigma$ in the sense that $\Phi(L) \equiv 0$, implies
that $\Phi$ is projective on the x-domain. When $Id=\delta_{x}(y)$, we have that Id has support on the diagonal 
$x=y$, that is a local operator. Given Id has point support and x fixed, we have $d I \notin G$, when $x \neq y$. 
In particular, when the domain is G, $u \neq v$ implies Id has no support, when u varies. This means that
an approximation property is only possible, by including the diagonal. Further, 
$\Phi \equiv 0$ on the diagonal (spirals). Example: $ u \rightarrow v$ as continuation,  the diagonal gives trivial continuations. 
Orthogonal continuations can be given as regular.

Consider $d V \rightarrow d V \times d V \rightarrow d U \times d U^{\bot}$. Starting from a strict
condition $d U^{\bot} / d U \rightarrow 0$, we get $ d U^{\bot}(V,W) \mid_{W=V}  \rightarrow 0$, when $d U \rightarrow d V$.
Example: consider $d U^{\bot}$ as a measure in $\{ d U \}^{\bot}$, we have $d U^{\bot} / d V \rightarrow 0$, that is a strictly weaker continuation,
in particular $d V$ and $d U^{\bot}$ are not of the same type. When V is Lagrange, we assume $\Delta V=0$, that is $(d U,d U^{\bot})$ corresponding to
real type relative d V. The condition $d U^{\bot} / d V \rightarrow 0$ regularly, implies that $(U,U^{\bot})$
has a trivial polar set relative (V,V). Example: assume $d V / d I = \rho \rightarrow 1$ regularly.
Further, $d V \rightarrow d U$ according to $\rho \rightarrow 1/\rho$, a sufficient condition for $d U^{\bot}/ d I \rightarrow 0$, is given  by
$d U^{\bot} / d U=\sigma \prec \prec \rho$, that is given $\rho \rightarrow \infty$, $\sigma (1/\rho^{2}) \rho \rightarrow 0$.
Further, $\rho \rightarrow 1$, as $t \rightarrow 0$, gives the condition for a diagonal, $\sigma / \rho \rightarrow 1$ regular.

Example: assume $\Phi$ projective, according to $d U = \rho d V$, with $\rho \neq 0$ on L.
Further, $\rho \rightarrow 0$ on L, with $\rho \neq 0 $ outside $\infty$. The corresponding condition $d U - d I=\rho (d V - d I)$,
gives that the type is preserved, the condition is dependent of topology.

\newtheorem{completion}[MP]{Topological completion}
\begin{completion}
When $\frac{d \mathcal{F}}{d x}$ is not bijective, consider $v$ so that $\frac{d \mathcal{F}}{d v}$ is bijective.
We get density through topological completion to R(V). 
\end{completion}

Note that $d \mathcal{F}(u) / d u \rightarrow 0$ regularly, when $u \rightarrow \infty$, that is $d v/ d u=0$ 
implies $1/u=0$, sufficient is $\mid d u / d v \mid \leq \mid u \mid$.

Example: Vf=0 iff Wf=f, that is start from invariant sets.
Assume $d I \in (K(G_{4}))_{H} \subseteq (G_{2})_{C}$. Given $f=T(\phi)$ we assume $f(u,u^{\bot}) \rightarrow f(0)$ in $C^{0}$
iff $f(u_{1},\ldots, u_{8}) \rightarrow f(0)$ in H. The principle for preservation of approximation property is
$\mathcal{H}_{1} \subset \mathcal{H}_{2}$ implies $G_{\mathcal{H}_{2}} \subset G_{\mathcal{H}_{1}}$.
Note that $(U,U^{\bot}) f \rightarrow f(u,u^{\bot})$ is not necessarily continuous. However $(U,U^{\bot})f = f(u_{1},\ldots,u_{8})$,
that is the domain for regular convergence in $G_{2}$ has the corresponding domain for the arguments in $G_{8}$. 
Isolated singularities in $G_{8}$ are according to $\{ g(v)=\delta_{j} g(v)=0 \}$ trivial, for $j=1,\ldots, 8$. 
The corresponding condition in $G_{2}$ is not continuous derivatives, but possible existence of w, so that $g(u,u^{\bot},w) \in C^{1}$. 

Assume $d I \notin G$, but $d U - d I \in G$, along invariant sets in the domain, 
we do not assume an approximation property, thus not all approximations are regular. 
d U - d I ($\simeq d U^{\bot}$) bijective corresponds to symplecticity. 
Given completion through $V(f)$, we have $\int f (d U- dI)=0$ implies $Vf=0$, 
that is trivial when V algebraic

N considered as a functional over d U, gives that $N \in T^{*}(bd X)$, where $X=R(d U)$.
Consider $N(d U) = d N_{1}(U)$, given $bd X$ strictly pseudo convex, we have that $\Sigma= \subset T^{*}(bd X)$
is symplectic.

\newtheorem{TS}[MP]{Semi-algebraic sets}
\begin{TS} 
Consider $d V \in G_{8}$ with $d V + d V^{L} = d U_{1} \in G_{2}$, with $d I \in G_{2}$,
Then we have $\{ f(v,v^{L}) < \lambda \}$ semi-algebraic, implies $\{ f(v,\cdot) < \lambda \}$ 
semi-algebraic.
\end{TS}

Assume $d I \in G_{2}$. For a geometrically convex closure of $If$ , we consider
$d V + d V^{L} \in G_{2}$. For holomorphic convex closure, we assume $d V + d V^{L}=d I$ iff $d W=0$ (analytic),
with $\mid f(w) \mid \leq \sup_{K} \mid  f \mid$, for a fixed compact K, defining a regular neighborhood of 
$d V + d V^{L}=d I$. Example: $\mid f(w) \mid \leq e^{\epsilon \mid v + v^{L} \mid}$ close to invariants. 
Given $d W$ projective, zero-lines L are preserved. Given
$d W$ algebraic according to =Pd I, for a polynomial P, we must have $L =\{ 0 \}$. 
Given d V projective, PdI algebraic corresponds to a strict condition.

Assume $d U^{\bot}=\rho d U$, where $\rho \rightarrow 1$ in $T_{ac}(H(G_{2}))$, that is the disk closure 
when $\rho \rightarrow 1$ in $C(G_{2}) \supset H(K(G_{4}))$, that is we assume $\rho \in H$ implies 
$\rho \neq 1$. In particular, $\{ G_{4} \quad \rho < \lambda \}$ semi-algebraic, but not necessarily relatively compact.

Weak type is defined through $<T,d U>(\phi)=<T,d I>(\phi)$, $\forall T \in (I)$, 
where (I) is a class of symbols, dependent of base $\phi$. Given $U \rightarrow {}^{t} U$ preserve 
invariant sets $\Sigma(U)$ over the base, that is $U,{}^{t} U$ of the same type, we can define 
invariant sets over the base. 

Consider, given $\{ d N = d I \}$ rectifiable, $d V \times d V^{L} \rightarrow d N \times d N^{\bot} \rightarrow d U \times d U^{\bot}$, 
so that $\{d V \times d V^{L}=d I \} \rightarrow \{ d N = d I \}$ absolute continuous and preserves sets of measure zero. 
We assume $(N,N^{\bot})$ defines a normal operator, so that for instance $N^{\bot} \phi=0$ implies 
$\phi$ has bounded support. Given $d V \rightarrow d {}^{t} V$ preserves type, we get decomposable sets $\{ d U=d I \}$ 
according to $(d V \otimes d I)(d I \otimes d V^{L})=d I \otimes d I$. Note that given $(d V,d V^{L})$ is dense in the domain
for $d U$ and projective as operators, we have $(d V \otimes d I)(d I \otimes d V^{L}) \simeq d V + d V^{L}$. Further, 
when $V^{L} \simeq {}^{t} V$ (harmonic), $UT(\phi)=VT(V \phi)$.
Invariant sets for strong topology, have corresponding invariant sets for weak topology. Convergence 
$d V \rightarrow d I$, implies $d {}^{t} V \rightarrow d I$, that is we have two-sided convergence. 
In particular, $d V \times d V^{L}$ preserve relatively compact sub-level surfaces. Conversely, 
the mapping $\{ d U^{\bot} =0 \} \rightarrow \{d N^{\bot} = 0 \}$ is absolute continuous and preserve 
sets of measure zero. Thus, $N(U,U^{\bot})$ is at the boundary given by N(U), that is not a spiral axis.

\subsubsection{Order of reducible movements}
\newtheorem{ordning}[MP]{Order for constant surfaces}
\begin{ordning}
Assume $A_{1}A_{2}A_{3}=I$ under the condition 
$A_{1}=A_{2}^{-1}$, then we have order 1, given all three $=I$, we have order 3 (maximal order)
\end{ordning}

Note $f(u_{1},u_{2})=const$ in C iff $f(u_{1},\ldots,u_{8})=const$ in H, that is the order is dependent of topology.
Consider $M(x/y,y/x)$, with dense domain and M algebraic, the condition $x/y \rightarrow 0$ or $y/x \rightarrow 0$ in $\infty$
implies not on the reflection axes, that is M is assumed algebraic outside invariant sets.
Irreducibels: $(A-I)(B-I)=I$ iff $AB=A+B$, that is irreducibels give decomposable 
constant surfaces. 

AB=I with $A \bot B$ non-linear, implies A=B=I that is of order 2.
Assume AB=I, given $A=1/B$ single valued in the plane, then A can be given linearly. Example: choose B monotonous 
and minimal, according to $d B \neq 0$.
Example: $d U^{-1} \notin (G_{2})_{H}$ and UV=I, with $d U,d V \in (G_{2})_{H}$ implies $U = V = I$, that is irreducibels
relative $(G_{2})_{H}$.

Example: $d U = \rho d I$, with $\rho \rightarrow c$ in $\infty$, $d W=\mu d I$ and $d U=\sigma d W$,
with $\sigma \rightarrow 0$ in $\infty$, we then have $\sigma=\rho / \mu$, sufficient is $\mu \rightarrow \infty$ in $\infty$.

\subsubsection{Continuation according to TP}

Continuation according to TP can be seen as a linear problem, given independence of $y'$.
We have independence of $y'$ (\cite{Lie91}), given
$\eta_{x}=0$, $(\eta_{y} - \xi_{x})=0$,$\xi_{y}=0$, for instance $(\xi,\eta)$ is given by 
$(1,1),(x,y)$.
Assume $d I \rightarrow d N \rightarrow d U=\rho d N$, where $d N \rightarrow d U$ continuous, 
in particular $d N \rightarrow d I$ implies $d N^{\bot} \rightarrow 0$,
thus N has order 2, that is translation, rotation.
We further assume $d U^{\bot}$ dense in $R(d U)^{\bot}$, in the sense that $\forall g \in \mathcal{H}$, we have existence of 
$f \in R(N)$, so that $(U + U^{\bot}) f=g$.

\newtheorem{reduktion}[MP]{Reduction to a normal model}
\begin{reduktion}
Consider $f_{L}(u,v)$, where $u,v$ are one-parameter, for instance $f_{L}$ is considered over a normal model. 
Given $\{ f_{L} \}$ has an approximation property, with $u \rightarrow 0$,$V \rightarrow 0$ (or $\delta f/\delta v \rightarrow 0$), we have that $\{ f \}$, 
a continuation according to TP, has an approximation property. In particular, $I f_{L}=f_{L} I$ corresponds to $f_{L}$ is topologically algebraic in $(u,v)$.
Assume $f=f_{L} + R$, where R corresponds to multivalentness for the representation of f, consider for instance $\{ R < \lambda \} \subset \subset X $,
a semi-algebraic set for R and $X \subset C^{\infty}$.
Given a compact neighborhood of this set, we can assume continuation to $C^{\infty}$,
that is $f-f_{L}$ can be continued to $C^{\infty}$.
\end{reduktion}

When $(u,v)$ convex and dense in the domain, we assume for a normal model, an approximation property 
$(u,v) \rightarrow 0$ regular, that is $f_{L} = \int f(u,v) dI(u,v) \rightarrow f(0)$ regularly, where we assume $u=v$ implies 
$v=0$. 

In general, we have $\int_{X} (d f) \omega = \int_{bd X} f \omega - \int_{X} f d \omega$. 
Topological degree gives the order for G that preserve bd X. Using the logarithm mapping, we have that lines are mapped
on to closed contours. Absence of entire lines in R(U=I), correspond to absence of closed curves (non-trivial), 
that preserve $d U=d I$.

\bibliographystyle{amsplain}
\bibliography{sos}

\providecommand{\bysame}{\leavevmode\hbox to3em{\hrulefill}\thinspace}
\providecommand{\MR}{\relax\ifhmode\unskip\space\fi MR }
\providecommand{\MRhref}[2]{%
  \href{http://www.ams.org/mathscinet-getitem?mr=#1}{#2}
}
\providecommand{\href}[2]{#2}
\begin{thebibliography}{10}

\bibitem{AhlforsSario60}
L.~Sario~L. Ahlfors, \emph{Riemann Surfaces}, Princeton University Press, 1960.

\bibitem{Ramis92}
J-P.~Ramis B.~Malgrange, \emph{Fonctions Multisommables}, Annales de l'institut
  Fourier \textbf{t. 42} (1992).

\bibitem{Bourbaki89}
N.~Bourbaki, \emph{General Topology, ch. 1-4}, Hermann, Paris, 1989.

\bibitem{Cartan41}
E.~Cartan, \emph{La Notion d'Orientation dans les Diff\'erentes
  G\'eom\'etries}, Bull. de la SMF \textbf{t. 69} (1941).

\bibitem{Collingwood66}
E.F. Collingwood and A.~J. Lohwater, \emph{The Theory of Cluster Sets},
  Cambridge Tracts in Math. and Math. Phys., no.56, 1966.

\bibitem{Dahn13}
T.~Dahn, \emph{Some Remarks on Tr{\`e}ves' Conjecture}, ArXiv (2013).

\bibitem{Dahn15}
\bysame, \emph{On Partially Hypoelliptic Operators, part i,ii}, ArXiv (2015).

\bibitem{Dahn18}
\bysame, \emph{Some Remarks on Schr{\"o}dinger Operators.}, ArXiv (2018).

\bibitem{Dahn2022}
\bysame, \emph{On a Weak Lie Type for Vectors}, ArXiv (2022).

\bibitem{BdMv71}
L.~Boutet de~Monvel, \emph{Utilisation des Op\`erateurs Pseudo-diff\`erentiels
  et des Op\`erateurs Trace dans Certains Probl\`emes aux Limites (suite)},
  SEDP \textbf{exp. no. 3} (1971).

\bibitem{BdMv85}
\bysame, \emph{Syst\`emes Presque Elliptiques: une Autre D\'emonstration de la
  Formule de l' Index}, Asterisque \textbf{131} (1985).

\bibitem{Riesz56}
B.~Sz.-Nagy {F. Riesz}, \emph{Functional Analysis.}, Dover Publications Inc.,
  1956.

\bibitem{Garding60}
L.~G{\aa}rding, \emph{Vecteurs Analytiques dans les Repr{\'e}sentations des
  Groupes de Lie}, Bulletin de la S.M.F \textbf{t.88} (1960).

\bibitem{Garding77}
L.~Gårding, \emph{Sharp Fronts of Paired Oscillatory Integrals}, RIMS, Kyoto
  \textbf{12} (1977).

\bibitem{Leray34}
J-Leray J.~Schauder, \emph{Topologie et \'Equation Fonctionelles}, Annales sci.
  de l'E.N.S. \textbf{3e s\'erie, t. 51} (1934).

\bibitem{Malgrange71}
B.~Malgrange, \emph{Ch. i. Forme Infinit\`esimale des \`Equations de Lie},
  Cours de l'institut Fourier \textbf{t. 8} (1971).

\bibitem{Nishino62}
T.~Nishino, \emph{Sur les Familles de Surfaces Analytiques}, J. Math. Kyoto
  Univ. \textbf{1-3} (1962).

\bibitem{Riesz20}
M.~Riesz, \emph{\"Uber die Randwerte einer Analytischen Funktion}, 4e Congres
  des Math. Scand. Stockholm (1920).

\bibitem{Riesz27}
\bysame, \emph{Sur les Maxima des Formes Bilin\'eaires et sur les Fonctionelles
  Lin\'eaires}, Acta Math. \textbf{49} (1927).

\bibitem{Lie91}
S.~Lie~G. Scheffers, \emph{Vorlesungen {\"u}ber Differentialgleichungen mit
  bekannten infinitesimalen Transformationen.}, Teubner Leipzig, 1891.

\bibitem{Schwartz57}
L.~Schwartz, \emph{Th{\'e}orie des Distributions {\`a} Valeurs Vectorielles,I},
  Ann. de l'institut de Fourier (1957).

\bibitem{Treves67}
F.~Treves, \emph{Topological Vector Spaces, Distributions and Kernels.},
  Academic Press, 1967.

\end{thebibliography}

\end{document}